\documentclass[twocolumn]{IEEEtran}
\IEEEoverridecommandlockouts
\usepackage{color}
\usepackage{booktabs, multirow, makecell}
\usepackage{bm, comment}
\usepackage[lofdepth,lotdepth]{subfig}
\usepackage{amsmath,amsfonts,amssymb}
\usepackage{mathabx,fixmath,algorithm, algorithmic}
\usepackage{pstool}
\usepackage{psfrag}
\usepackage{hyperref}
\usepackage{cite}

\def\x{\mathbold{x}}
\def\A{\mathbold{A}}

\def\L{\mathbold{L}}

\def\blambda{\mathbold{\lambda}}
\def\bmu{\mathbold{\mu}}
\def\bxi{\mathbold{\xi}}
\def\cf{f^{\star}}

\def\Q{\mathbf{Q}}

\def\v{\mathbold{v}}
\def\c{\mathbold{c}}
\def\p{\mathbold{p}}

\def\b{\mathbold{b}}

\def\y{\mathbold{y}}

\newcommand{\im}{\mathrm{im}}
\newcommand{\nul}{\mathrm{null}}
\newcommand{\reals}{\mathbb{R}}

\def\transp{\mathsf{T}}

\def\qed{\hfill $\blacksquare$ }
\DeclareMathOperator*{\argmin}{\textrm{argmin}}

\usepackage{theorem}

\newtheorem{theorem}{Theorem}
\newtheorem{corollary}{Corollary}
\newtheorem{proposition}{Proposition}
\newtheorem{remark}{Remark}
\newtheorem{assumption}{Assumption}
\newtheorem{claim}{Claim}
\newtheorem{example}{\bf Example}

\newcommand{\cref}[1]{\eqref{#1}}

\usepackage{lipsum}
\usepackage{amsfonts}
\usepackage{graphicx}
\usepackage{epstopdf}
\usepackage{algorithmic}
\ifpdf
  \DeclareGraphicsExtensions{.eps,.pdf,.png,.jpg}
\else
  \DeclareGraphicsExtensions{.eps}
\fi

\begin{document}
\title{Dual Prediction-Correction Methods for Linearly Constrained Time-Varying Convex Programs} 
\author{Andrea Simonetto

\thanks{

A preliminary version of this work with limited results and no proofs has been submitted to the American Control Conference 2018 as~\cite{ACC}.

Andrea Simonetto is with the Optimization and Control Group of IBM Research Ireland, Dublin, Ireland. Email: andrea.simonetto@ibm.com . }
}

\maketitle

\begin{abstract}
Devising efficient algorithms to solve continuously-varying strongly convex optimization programs is key in many applications, from control systems to signal processing and machine learning. In this context, solving means to find and track the optimizer trajectory of the continuously-varying convex optimization program. 
Recently, a novel prediction-correction methodology has been put forward to set up iterative algorithms that sample the continuously-varying optimization program at discrete time steps and perform a limited amount of computations to correct their approximate optimizer with the new sampled problem and predict how the optimizer will change at the next time step. Prediction-correction algorithms have been shown to outperform more classical strategies, i.e., correction-only methods. Typically, prediction-correction methods have asymptotical tracking errors of the order of $h^2$, where $h$ is the sampling period, whereas classical strategies have order of $h$. Up to now,  Prediction-correction algorithms have been developed in the primal space, both for unconstrained and simply constrained convex programs. In this paper, we show how to tackle linearly constrained continuously-varying problem by prediction-correction in the dual space and we prove similar asymptotical error bounds as their primal versions. 
\end{abstract}



\section{Introduction}\label{sec:intro}

Continuously varying optimization programs have appeared as a natural extension of time-invariant ones when the cost function, the constraints, or both, depend on a time parameter and change continuously in time. This setting captures relevant control, signal processing, and machine learning problems (see e.g.,~\cite{ACC} for a broad overview).

We focus here on linearly constrained time-varying \emph{convex} programs of the form
\begin{equation}\label{eq.tvar}
\x^*(t) : = \argmin_{\x \in \mathbb{R}^n} \, f(\x; t), \quad \textrm{subject to: } \A\x = \b,
\end{equation}
where $t\in \reals_+$ is non-negative, continuous, and it is used to index time;  $f: \mathbb{R}^n\times \mathbb{R}_{+} \to \mathbb{R}$ is a \emph{smooth strongly convex function} uniformly in time; $\A \in \mathbb{R}^{p \times n}$ and $\b \in \mathbb{R}^m$ are a real-valued matrix and vector that represent the linear equality constraints. The goal is to find (and track) the solution $\x^*(t)$ of~\cref{eq.tvar} for each time $t$ -- hereafter referred to as the  optimal solution \emph{trajectory}. 

Problem \cref{eq.tvar} might be solved in a centralized setting based on a continuous time platform \cite{Fazlyab2016}; however, here we focus on a discrete time setting. The reason for this choice is motivated by the widespread use of digital computing units, such as control units (actuators) and digital sensors. In this context, we envision that our optimization problem will change in response to measurements taken at discrete time steps and its solution could provide control actions to be implemented on digital control units, similarly to~\cite{Zavala2010}. In addition, we also envision our methods to be implemented on networks of communicating and computing nodes. In this latter scenario, at each time step the nodes have to send messages among each other before and during the computations. Then, continuous-time settings may be less appropriate, especially when high communication latencies may be expected. 

Therefore, we use sampling arguments to reinterpret~\cref{eq.tvar}  as a sequence of time-invariant problems. In particular, upon sampling the objective functions $ f(\x; t)$ at time instants $t_{k}$,  $k=0,1,2,\dots$, where the sampling period $h := t_k - t_{k-1}$ can be chosen arbitrarily small, one can solve the sequence of time-invariant problems
\begin{equation}\label{eq.tinv}
\x^*(t_k) : = \argmin_{\x \in \mathbb{R}^n} \, f(\x; t_k), \quad \textrm{subject to: } \A\x = \b.
\end{equation}
By decreasing $h$, an arbitrary accuracy may be achieved when approximating problem~\cref{eq.tvar} with \cref{eq.tinv}. However, solving \cref{eq.tinv} for each sampling time $t_k$ may not be computationally affordable in many application domains, even for moderate-size problems. 

Focusing on unconstrained or simply constrained optimization problems, a series of works among which~\cite{Paper2, Paper3} developed prediction-correction methods to find and track the solution trajectory $\x^*(t)$ up to a bounded asymptotical error, in the \emph{primal space}. This methodology arises from non-stationary optimization~\cite{Polyak1987}, parametric programming~\cite{Robinson1980,Zavala2010,Dontchev2013,Kungurtsev2017}, and continuation methods in numerical mathematics~\cite{Allgower1990}. 

This paper extends the current state-of-the-art methods~\cite{Zavala2010,Fazlyab2016,Paper3} by offering the following contributions. 

First, we develop prediction-correction methods to track the solutions of the time-varying \emph{linearly constrained} problems~\cref{eq.tvar} by leveraging a dual ascent technique. To the author's knowledge, this is the first work that proposes prediction-correction methods in the \emph{dual space}. In~\cite{Jakubiec2013, Simonetto20XX}, the authors have developed dual ascent methods for similar problems, but they are correction-only methods and -- as we prove here -- they have worse tracking capabilities than dual prediction-correction methods. 


Second, our algorithm can handle a \emph{rank deficient} matrix $\A$, which is a situation ubiquitous in distributed optimization problems. This therefore opens the way to distributed algorithms based on dual decomposition, which have many applications. This was not considered in previous efforts, i.e., in the continuous-time platform of~\cite{Fazlyab2016}\footnote{The work in~\cite{Fazlyab2016} differs from the work here, not only because they work in continuous time and they do not consider a rank deficient  $\A$. It differs also from the algorithmic perspective: they propose a continuous-time primal-dual algorithm, while here we focus on discrete-time dual ones. }. 

In this paper, we derive methods that are proved to track the solution trajectory $\x^*(t)$ up to an asymptotical error upper bound, which depends on the sampling period, on the properties of the cost function, and on the number of prediction and correction steps we use, and on the spectral properties of $\A$.
%
%
%
With the aid of numerical simulations, we are able to showcase further the performance of the proposed methods and their comparison with the correction-only strategies. In particular, the proposed algorithms outperform the correction-only ones in asymptotic error bounds and they appear also better when computational considerations are taken into account, in most cases.



{\bf Organization.} In Section~\ref{sec:cases}, we introduce the basic assumptions for the linear system $\A\x = \b$. Section~\ref{sec:timeinvariant} covers the required background on time-invariant and correction-only methods in the dual domain. We present our algorithm in Section~\ref{sec:pc}, while convergence analysis is discussed in Section~\ref{sec:conv}. The main result of this paper is presented in Theorem~\ref{theo.main}. Section~\ref{sec:distr} studies distributed optimization problems. Numerical simulations are presented in Section~\ref{sec:num}, and we conclude in Section~\ref{sec:conc}. 


{\bf Notation.} Vectors are written as $\x\in\reals^n$ and matrices as $\A\in\reals^{p\times n}$. We use $\|\cdot\|$ to denote the Euclidean norm in the vector space, and the respective induced norms for matrices and tensors. The image (i.e., the column space) and the nullspace of matrix $\A$ are indicated as $\im(\A)$ and $\nul(\A)$ respectively. The gradient of the function $f(\x; t)$ with respect to $\x$ at the point $(\x,t)$ is denoted as $\nabla_{\x} f(\x; t) \in \reals^n$, the partial derivative of the same function with respect to (w.r.t.) $t$ at $(\x,t)$ is written as $\nabla_t f(\x; t)\in \reals$. Similarly, the notation $\nabla_{\x\x} f(\x; t) \in \reals^{n\times n}$ denotes the Hessian of $f(\x;t)$ w.r.t. $\x$ at $(\x,t)$, whereas $\nabla_{t\x} f(\x; t) \in \reals^{n}$ denotes the partial derivative of the gradient of $f(\x;t)$ w.r.t. the time $t$ at $(\x,t)$, i.e. the mixed first-order partial derivative vector of the objective. The tensor $\nabla_{\x\x\x} f(\x; t) \in \reals^{n\times n\times n}$ indicates the third derivative of $f(\x;t)$ w.r.t. $\x$ at $(\x,t)$, the matrix $\nabla_{\x t\x} f(\x; t) = \nabla_{t\x\x} f(\x; t)  \in \reals^{n\times n}$ indicates the time derivative of the Hessian of $f(\x;t)$ w.r.t. the time $t$ at $(\x,t)$, the vector $\nabla_{t t\x} f(\x; t) \in \reals^{n}$ indicates the second derivative in time of the gradient of $f(\x;t)$ w.r.t. the time $t$ at $(\x,t)$.


\section{Assumptions for $\A\x = \b$}\label{sec:cases}

We assume that $\b \in \im(\A)$, so that the optimization problems~\cref{eq.tvar} (or equivalently~\cref{eq.tinv}) has a solution for each time $t$. We do not assume that the matrix $\A\in\mathbb{R}^{p\times n}$ is full row rank, so it can be rank deficient. 

The singular values of $\A$ are ordered as $\sigma_{\max}:=\sigma_p \geq \sigma_{p-1} \geq \sigma_{\min} > \sigma_j = \dots = \sigma_1 = 0$, where $\sigma_{\min}$ is the minimum positive singular value. We call $\kappa_\A = \sigma_{\max}/\sigma_{\min}$.

Since $\b \in \im(\A)$, one could eliminate the redundant rows in $\A$ (if rank deficient) and construct a full row rank matrix. However, in some cases it is desirable to keep a rank deficient $\A$, since it encodes more linear constraints. This is the case, e.g., in distributed optimization where $\A$ describes the communication links that are present. The more links means (in general) faster convergence to the desired solution.

\section{Time-invariant and correction-only dual ascent}\label{sec:timeinvariant}

We start by two (known) properties of the primal and dual variables at optimality, for the time-invariant problem~\cref{eq.tinv}.  

\begin{proposition}\label{prop:properties}
Let function $f(\cdot; t_k): \mathbb{R}^n \to \mathbb{R}$ be strongly convex with constant $m$ and strongly smooth with constant $L$. Then the primal optimizer of~\cref{eq.tinv}, i.e., $\x^*(t_k)$, is unique. 

If $\A$ is full row rank, the dual optimizer of~\cref{eq.tinv}, i.e., $\blambda^*(t_k)$, is also unique and $\blambda^*(t_k) \in \im(\A)$. 

If $\A$ is rank deficient, there exists a unique dual optimizer of~\cref{eq.tinv} $\blambda^*(t_k)$ for which $\blambda^*(t_k) \in \im(\A)$.
\end{proposition}

\begin{proof} Given in Appendix~\ref{appendix:uno}. \end{proof}

Proposition~\ref{prop:properties} sets the frame for the results of this paper. If $\A$ is full row rank, the primal-dual optimizers are unique and the dual optimizer lies in the image of $\A$. If $\A$ is rank deficient,  the primal optimizer is unique, while the dual is not unique but we will be interested in finding the unique dual optimizer that lies in the image of $\A$. By restricting the search space to the image of $\A$, we will be able to overcome the rank deficiency of $\A$ in the proofs, without losing optimality.   

Consider now the following iterative algorithm to solve~\cref{eq.tinv}, known as dual ascent. 

\begin{enumerate}
\item Pick $(\x_0, \blambda_0)$; Set $i = 0$; Pick a stepsize $\alpha>0$.
\item Iterate: 
\begin{subequations}\label{dual.ascent}
\begin{align}\label{dual.ascent1}
\x_{i+1} &= \argmin_{\x \in \mathbb{R}^n}\{f(\x; t_k) + \blambda_{i}^\transp \A\x \},\\
\blambda_{i+1} &= \blambda_{i} + \alpha (\A\x_{i+1} - \b).
\end{align}
\end{subequations}
\end{enumerate}

We have the following result. 


\begin{theorem}[Time-invariant dual ascent convergence]\label{theo.dualascent}
Fix the time $t_k$. Let function $f(\cdot;t_k): \mathbb{R}^n \to \mathbb{R}$ be strongly convex with constant $m$ and strongly smooth with constant $L$. Select $\x_0$ arbitrarily, but $\blambda_0 \in \im(\A)$. 
%
%
%
%
%
Let the stepsize $\alpha$ be chosen as $\alpha < 2 m/\sigma_{\max}^2$. Then, the sequence $\{(\x_i, \blambda_i)\}_{i \in \mathbb{N}}$ generated by recursively applying~\cref{dual.ascent} converges to the unique primal-dual optimizer of~\cref{eq.tinv} $(\x^*(t_k), \blambda^*(t_k) \in \im(\A))$. In particular, $\{\blambda_i\}_{i\in\mathbb{N}}$ converges Q-linearly to $\blambda^*(t_k) \in \im(\A)$ as 
\begin{equation}
\|\blambda_{i+1} - \blambda^*(t_k)\| \leq \varrho \|\blambda_{i} - \blambda^*(t_k)\| \leq \varrho^{i+1} \|\blambda_{0} - \blambda^*(t_k)\|, \label{claim1}
\end{equation}
while $\{\x_i\}_{i\in\mathbb{N}}$ converges R-linearly as 
\begin{equation}
\|\x_{i+1} - \x^*(t_k)\| \leq \frac{\sigma_{\max}}{m} \|\blambda_{i} - \blambda^*(t_k)\|,\label{claim2}
\end{equation}
where the contraction factor $\varrho < 1$ is defined as $\varrho = \max\{|1-\alpha \sigma_{\max}^2/m|, |1-\alpha \sigma_{\min}^2/L|\}$.
\end{theorem}

\begin{proof} Given in Appendix~\ref{appendix:due}. \end{proof}


Theorem~\ref{theo.dualascent} says that the \emph{time-invariant} iteration~\cref{dual.ascent} converges to the primal-dual optimizer of the \emph{time-invariant} optimization problem~\cref{eq.tinv}. Furthermore, the rate is linear. 

In~\cite{Jakubiec2013,Simonetto20XX}, the authors extend the previous results to a running version (or with the nomenclature here, a correction-only version) of dual ascent. By running we mean an algorithm that adjust the problem on-line while the algorithm is running. In this context, consider the time-varying problem~\cref{eq.tvar} and the running version of the iterations~\cref{dual.ascent} defined by sampling problem~\cref{eq.tvar} at discrete sampling times, as follows:

\begin{enumerate}
\item Pick $(\x_0, \blambda_0)$; Set $k = 0$; Pick a stepsize $\alpha>0$.
\item Iterate: 
\begin{subequations}\label{dual.ascent.tv}
\begin{align}\label{dual.ascent1.tv}
\x_{k+1} &= \argmin_{\x \in \mathbb{R}^n}\{f(\x; t_k) + \blambda_{k}^\transp \A\x \},\\
\blambda_{k+1} &= \blambda_{k} + \alpha (\A\x_{k+1} - \b).
\end{align}
\end{subequations}
\end{enumerate}

As one can see, this running version of~\cref{dual.ascent} considers functions that change at the same time as the updates are computed (i.e., there is only one time variable $k$). 
 
Let the following assumptions hold. 

\begin{assumption}\label{as.1}
Let time-varying function $f: \mathbb{R}^n\times \mathbb{R}_{+} \to \mathbb{R}$ be strongly convex with constant $m$ and strongly smooth with constant $L$, uniformly in time (i.e., for each time $t\geq 0$). Define the condition number of $f$ as $\kappa_f:=L/m$, uniformly in time.
\end{assumption}

\begin{assumption}\label{as.2}
Let the distance between optimizers of Problem~\cref{eq.tvar} at two subsequent sampling time $t_k$ and $t_{k-1}$, i.e., $(\x^*(t_k)$, $\blambda^*(t_k)\in \im(\A))$ and $(\x^*(t_{k-1})$, $\blambda^*(t_{k-1})\in \im(\A))$, be upper bounded for each $k>0$ as,
\begin{equation}
\max\{\|\x^*(t_k)-\x^*(t_{k-1})\|, \|\blambda(t_k)^*-\blambda^*(t_{k-1})\|\} \leq K. 
\end{equation}
\end{assumption}
 

Then the following result is in place. 
\begin{theorem}[Running dual ascent convergence]\label{theo.dualascent.tv}
Under Assumptions~\ref{as.1}-\ref{as.2}, consider the running iterations~\cref{dual.ascent.tv}. Select $\x_0$ arbitrarily, but $\blambda_0 \in \im(\A)$.
Let the stepsize $\alpha$ be chosen as $\alpha < 2 m/\sigma_{\max}^2$. Then, the sequence $\{(\x_k, \blambda_k)\}_{k\in\mathbb{N}}$ generated by recursively applying~\cref{dual.ascent.tv} converges to the unique primal-dual trajectory of~\cref{eq.tvar}, $(\x^*(t_k), \blambda^*(t_k) \in \im(\A))$, up to a constant error bound linearly as 
\begin{align}
\|\blambda_{k+1} - \blambda^*(t_{k+1})\| &\leq \varrho (\|\blambda_{k} - \blambda^*(t_k)\| + K)  \nonumber \\ & \leq \varrho^{k+1} \|\blambda_{0} - \blambda^*(t_0)\| + \frac{\varrho K}{1-\varrho}, \label{claim1.tv}\\
\|\x_{k+1} - \x^*(t_{k+1})\| &\leq \frac{\sigma_{\max}}{m} (\|\blambda_{k} - \blambda^*(t_k)\|+K),\label{claim2.tv}
\end{align}
where the contraction factor $\varrho < 1$ is defined as $\varrho = \max\{|1-\alpha \sigma_{\max}^2/m|, |1-\alpha \sigma_{\min}^2/L|\}$.
\end{theorem}

\begin{proof} The proof is given for example in~\cite{Simonetto20XX}, and it is based on the results of Theorem~\ref{theo.dualascent} and the triangle inequality. In particular, for each time $t_{k+1}$ one can write
\begin{multline}
\|\blambda_{k+1} - \blambda^*(t_{k+1})\| \leq  \varrho(\|\blambda_{k} - \blambda^*(t_{k+1})\|) \leq \\ \leq  \varrho ( \|\blambda_{k} - \blambda^*(t_k)\| + K),
\end{multline}
which is yield by directly applying the time-invariant results and the triangle inequality, and by leveraging Assumption~\ref{as.2} on the variability of the optimizers. 
\end{proof}

Theorem~\ref{theo.dualascent.tv} is a generalization of Theorem~\ref{theo.dualascent} for cases in which the cost function changes continuously in time. The convergence result is similar to those of Theorem~\ref{theo.dualascent} but is achieved up to a constant error bound, which is due to the drifting of the primal-dual optimal pair. In the limit, one obtains the bounds, 
\begin{subequations}
\begin{align}
\limsup_{k\to\infty}\|\blambda_{k} - \blambda^*(t_{k})\| & = \frac{\varrho K}{1-\varrho}, \\
\limsup_{k\to\infty}\|\x_{k} - \x^*(t_{k})\| &= \frac{\sigma_{\max}}{m} \left(\frac{\varrho K}{1-\varrho} + K\right),\label{claim2.tv3}
\end{align}
\end{subequations}
and when $K = 0$, i.e., we are back in the time-invariant scenario, one re-obtains exact convergence. 


\begin{remark}
The result presented in Theorem~\ref{theo.dualascent.tv} can be readily related back to those in~\cite{Jakubiec2013}. There, Assumption~\ref{as.2} is substituted with an assumption on the primal optimizers variation and their gradients. Given the optimal conditions for~\cref{eq.tvar}, one can always translate the latter in the former as
\begin{multline}
\|\nabla_{\x} f(\x^*(t_{k}); t_k) - \nabla_{\x} f(\x^*(t_{k-1}); t_{k-1})\|\leq \\ \leq\sigma_{\max} \|\blambda(t_k)^*-\blambda^*(t_{k-1})\| \leq \sigma_{\max} K.
\end{multline}
\end{remark}


\section{Prediction-correction methodology}\label{sec:pc}

The running dual ascent~\cref{dual.ascent.tv} is agnostic of variations of the cost function, in the sense that it only reacts to variations of the cost but it does not attempt at predicting how the function changes depending either on past data, or on the knowledge of the time derivatives of the function. Recently, e.g.,~\cite{Paper2,Paper3}, a prediction-correction methodology in discrete time has been put forward to increase the accuracy of running (i.e., correction-only) methods by predicting how the cost function changes in time. The aforementioned works stay in the primal space, while here we will extend them to the dual space. 

\subsection{Prediction step}

To develop the prediction step, we consider the optimality conditions for~\cref{dual.ascent.tv} at time $t_{k+1}$, 
\begin{equation}\label{kkt1}
\nabla_{\x} f(\x^*(t_{k+1}); t_{k+1}) + \A^\transp \blambda^*(t_{k+1}) = {\bf 0}, \, \A \x^*(t_{k+1})  = \b.
\end{equation}
At time $t_k$, one cannot solve~\cref{kkt1} to determine the next primal-dual pair. What one can do is to approximate~\cref{kkt1} with the knowledge one has at $t_k$ via a backward Taylor expansion as,
\begin{subequations} \label{taylor}
\begin{align}\label{kkt2}
\nabla_{\x} f(\x^*(t_{k+1}); t_{k+1}) + \A^\transp \blambda^*(t_{k+1}) & \approx  \nabla_{\x} f(\x^*(t_{k}); t_{k}) + \nonumber \\ & \hspace*{-5cm}+ \nabla_{\x\x} f(\x^*(t_{k}); t_{k}) \delta \x + \nabla_{t\x} f(\x^*(t_{k}); t_{k}) h + \nonumber \\  & \hspace{-2cm} +\A^\transp (\blambda^*(t_{k}) + \delta \blambda) =  {\bf 0} \\ & \hspace{-4cm}\A \x^*(t_{k+1}) =   \A (\x^*(t_{k}) +\delta \x)  = \b.
\end{align}
\end{subequations}
If then, one is provided with the primal-dual optimizers at time $t_k$, one can approximate (or predict) the next primal-dual optimal pair by solving~\cref{taylor} for $\delta\x$ and $\delta \blambda$. That is, one has to solve the following quadratic program
\begin{subequations}\label{qp}
\begin{align}
\min_{\delta\x \in\mathbb{R}^n} & \,  \frac{1}{2}\delta\x ^\transp \nabla_{\x\x} f(\x^*(t_{k}); t_{k})\delta\x \!+\!  h\,\nabla_{t\x} f(\x^*(t_{k}); t_{k})^\transp \delta\x, \\ & \textrm{subject to } \A\delta\x = {\bf 0}.
\end{align}
\end{subequations}
We use this reasoning to develop our prediction step. 

Let $(\x_k, \blambda_k)$ be an approximate primal-dual optimal pair available at time $t_k$. In the prediction step we solve the quadratic problem:
\begin{subequations}\label{qp1}
\begin{align}
\min_{\delta\x \in\mathbb{R}^n} &\,\frac{1}{2}\delta\x ^\transp \nabla_{\x\x} f(\x_k; t_{k})\delta\x \!+\!  h\,\nabla_{t\x} f(\x_k; t_{k})^\transp \delta\x, \\ &\textrm{subject to } \A\delta\x = {\bf 0},
\end{align}
\end{subequations}
and we set the predicted pair to $\x_{k+1|k} = \x_k + \delta \x$, $\blambda_{k+1|k} = \blambda_k + \delta \blambda$. 

To solve~\cref{qp1}, various techniques can be applied. If one has access to the full instance, one can find the unique $(\x_k, \blambda_k\in\im(\A))$ by solving~\cref{qp1} at optimality\footnote{This can be done e.g., by off-the-shelf solvers, or by a custom-made Newton's method which can also employ Krylov-subspace based solvers for the resulting linear system. Note that some of the computations could be made off-line since $\A$ is time-invariant.}.  

Since we would like to devise algorithms that can be implemented in a distributed way, we follow another approach, which is to set up a dual gradient method with the iterations:

\begin{enumerate}
\begin{subequations}\label{eq:dual.ascent.tv}
\item Pick $(\delta \x_0, \delta \blambda_0={\bf 0})$; Set $p = 0$; Pick a stepsize $\beta>0$ and a maximum number of iterations $P$.
\item Iterate till $p = P-1$:
\begin{align}
\delta \x_{p+1} &= \argmin_{\delta \x \in \mathbb{R}^n}\Big\{\frac{1}{2}\delta\x^\transp\nabla_{\x\x}f(\x_k; t_k) \delta\x + \nonumber \\ & \hspace{1.5cm}+ h\,\nabla_{t\x}f(\x_k; t_k)^\transp \delta\x+ \delta \blambda_{p}^\transp \A\delta \x \Big\},\label{dual.pred.opt}\\
\delta \blambda_{p+1}& = \delta \blambda_{p} + \beta \A\delta \x_{p+1}.
\end{align}
\item Set $\widehat{\x}_{k+1|k} = \x_k + \delta \x_{P}$, $\widehat{\blambda}_{k+1|k} = \blambda_k + \delta \blambda_P$
\end{subequations}
\end{enumerate}

This converges to the exact $\x_{k+1|k}$ and $\blambda_{k+1|k}$ as $P \to \infty$, due to Theorem~\ref{theo.dualascent}.


This last option (which determines the solution of~\cref{qp1} only approximately if $P$ stays finite -- and that is why we indicate the predicted variable with an hat) is to be preferred in \emph{distributed} settings (as we will see in Section~\ref{sec:distr}). Of course, to make this last option viable, the maximum number of iterations $P$ needs to be small enough, which will induce an extra error in computing the prediction step. 

For the sake of uniformity, from now on, we will indicate with $\widehat{\x}_{k+1|k}$ and $\widehat{\blambda}_{k+1|k}$ both the exact and approximate prediction: in fact, the exact prediction couple is equivalent to the approximate one when $P \to \infty$.

\subsection{Correction step}

At time $t_{k+1}$, when one is allowed to sample the new cost function $f(\cdot; t_{k+1})$, then a correction step can be performed starting from the (approximate or exact) predicted pair previously computed. The correction step is nothing else than one (or possible multiple) round(s) of the dual ascent iteration as
\begin{subequations}\label{dual.ascent.tv.correction}
\begin{enumerate}
\item Pick $(\v_0 = \widehat{\x}_{k+1|k},{\bxi_0}=\widehat{\blambda}_{k+1|k})$; Set $c = 0$; Pick a stepsize $\alpha>0$ and a maximum number of iterations $C$.
\item Iterate till $c = C-1$:
\begin{align}
\v_{c+1} &= \argmin_{\v \in \mathbb{R}^n}\{f(\v; t_{k+1}) + {\bxi}_{c}^\transp \A\v \}, \label{dual.corr.opt}\\
\bxi_{c+1} &= {\bxi}_{c} + \alpha (\A\v_{c+1} - \b).
\end{align}
\item Set ${\x}_{k+1} = \v_C$, ${\blambda}_{k+1} = \bxi_C$.
\end{enumerate}
\end{subequations}
 
\subsection{Complete algorithms}

In Algorithm~\ref{alg.approx}, we summarize the prediction-correction methodology for the approximate prediction. The algorithm is parametrized over the number of prediction and correction steps that it employs. 

\begin{algorithm}[t]
\caption{Approx. Dual Prediction-Correction (ADuPC)}
\label{alg.approx}
\textbf{Require:} Initial guess $(\x_0,\blambda_0\in \im(\A))$; stepsizes $\alpha, \beta>0$; number of prediction and correction steps $P, C$
\begin{algorithmic}[1]
\FOR{$k = 0,1,2,\ldots$} 
	\STATE \color{blue}// time $t_k$\color{black}
	\STATE \textbf{Prediction step}: Compute $\delta\x$ and $\delta\blambda$ by approximately solving the quadratic program~\cref{qp1} by using the iterations~\cref{eq:dual.ascent.tv} with $(\delta\x_0 = {\bf 0},\delta{\blambda_0}={\bf 0})$, stepsize $\beta$, and number of iterations $P$
	\STATE Set $\widehat{\x}_{k+1|k} = \x_k + \delta \x$, $\widehat{\blambda}_{k+1|k} = \blambda_k + \delta \blambda$ 
	\STATE \color{blue}// time $t_{k+1}$\color{black}
	\STATE Acquire the updated function $f(\cdot; t_{k+1})$
	\STATE \textbf{Correction step}: Compute $\x_{k+1}$ and $\blambda_{k+1}$ by using the iterations~\cref{dual.ascent.tv.correction} with $(\v_0 = \widehat{\x}_{k+1|k},{\bxi_0}=\widehat{\blambda}_{k+1|k})$, stepsize $\alpha$, and number of iterations $C$
\ENDFOR
\end{algorithmic}
\end{algorithm}

In the next section, we study the convergence of Algorithm~\ref{alg.approx} to a ball around the optimal primal-dual trajectory. The size of the error ball will depends on the sampling period and the number of prediction and correction steps, among other parameters. 


\section{Convergence analysis}\label{sec:conv}

To derive our convergence results, we need the following additional assumptions. 

\begin{assumption}\label{as.3}
The time derivative of the gradient of the cost function is uniformly upper bounded for all $\x \in \mathbb{R}^n$ as
\begin{equation*}
\|\nabla_{t\x} f(\x; t)\|\leq C_0, \quad \forall \x \in \mathbb{R}^n, t.
\end{equation*} 
\end{assumption}

\begin{assumption}\label{as.4}
The cost function has bounded third order derivatives with respect to $\x$ and $t$ as
\begin{align*}
\|\nabla_{\x\x\x} f(\x; t)\|\leq C_1, \,& \|\nabla_{\x t\x} f(\x; t)\|\leq C_2, \\  \|\nabla_{tt\x} f(\x; t)\|\leq C_3,& \quad \forall \x \in \mathbb{R}^n, t.
\end{align*}
\end{assumption}

Assumptions~\ref{as.3}-\ref{as.4} are common in the time-varying optimization domain when dealing with prediction-correction methods, see~\cite{Dontchev2013,Paper3,Fazlyab2016}.



Central to our analysis is the following novel implicit function theorem. 

\begin{theorem}[Implicit function theorem for Problem~\cref{eq.tvar}]\label{theo.implicit}
Consider the time-varying problem~\cref{eq.tvar}. Let Assumptions~\ref{as.1} and~\ref{as.3} hold. The primal-dual optimal trajectory $\{\x^*(t_{k}), \blambda^*(t_{k})\in\im(\A)\}$ is locally Lipschitz in time (i.e., for small enough sampling periods), and in particular,
\begin{subequations}\label{bounds.xl}
\begin{align}
\|\x^*(t_{k}) - \x^*(t_{k-1})\| &\leq \frac{\kappa_f \kappa_{\A}^2 + 1}{m}\,{C_0}\, h  = O(h), \label{bounds.xl.primal} \\
\|\blambda^*(t_{k}) - \blambda^*(t_{k-1})\| &\leq  \frac{\kappa_f \kappa_{\A}}{\sigma_{\min}}\, C_0\, h = O(h). \label{bounds.xl.dual}
\end{align}
\end{subequations}
In addition, if the bounds $C_1, C_2, C_3$ are all identically $0$, then the inequalities~\cref{bounds.xl} are valid globally (i.e., the trajectory is globally Lipschitz in time, i.e.,~\eqref{bounds.xl} are valid for all sampling periods). 
\end{theorem}
 
\begin{proof} Given in Appendix~\ref{appendix:tre}. \end{proof}

Theorem~\ref{theo.implicit} characterizes how the optimal primal-dual pair changes over time due to functional changes. In particular, Theorem~\ref{theo.implicit} implies that optimizers changes are Lipschitz continuous in time, for sufficiently small sampling periods. As we see, Theorem~\ref{theo.implicit} does not need Assumption~\ref{as.2}, which is substituted by the stronger Assumption~\ref{as.3}. In particular, one can see that Assumption~\ref{as.2} is automatically enforced, as follows.  

\begin{corollary}\label{co:dualascent_h}
Let Assumption~\ref{as.1} and~\ref{as.3} hold. Then Assumption~\ref{as.2} is automatically satisfied with 
\begin{equation}
K = \max\left\{\frac{\kappa_f \kappa_{\A}^2 + 1}{m}, \frac{\kappa_f \kappa_{\A}}{\sigma_{\min}}\right\}\, C_0\, h.
\end{equation}
In addition, the asymptotical error for the running dual ascent~\cref{dual.ascent.tv} is $O(h)$. 
\end{corollary}

Corollary says that the error bound of the running version of dual ascent is proportional to the sampling period $h$ whenever Assumptions~\ref{as.1} and \ref{as.3} hold.   

We are now ready to prove convergence of the approximate dual prediction-correction algorithm. 

\begin{theorem}[Convergence of Algorithm~\ref{alg.approx}] \label{theo.main}
Consider the time-varying problem~\cref{eq.tvar}. Let Assumptions~\ref{as.1} and~\ref{as.3} hold. Consider $P$ prediction steps and $C$ correction steps, while let the stepsizes for prediction $\beta$ and correction $\alpha$ be chosen such that $\beta< 2 m/\sigma_{\max}^2$,  $\alpha< 2m/\sigma_{\max}^2$. Define the contraction factors for prediction and correction as,
\begin{subequations} 
\begin{align}
\varrho_{\textrm{P}}&:=\max\{|1-\beta \sigma_{\max}^2/m|, |1-\beta \sigma_{\min}^2/L|\}, \\ \varrho_{\textrm{C}}&:=\max\{|1-\alpha \sigma_{\max}^2/m|, |1-\alpha \sigma_{\min}^2/L|\}.
\end{align}
\end{subequations}
Select the prediction and correction steps $P, C$ to verify the contraction property
\begin{equation}\label{cond:1}
\gamma_1:=\varrho_{\textrm{C}}^{C}(2\varrho_\textrm{P}^{P}+1)<1.
\end{equation}
There exists a constant $\gamma_2>0$, dependent on the problem parameters, such that if one chooses the sampling period $h$ as
\begin{equation}\label{cond:2}
h < (1-\gamma_1)/\gamma_2,
\end{equation}
%
%
(so that $\tau(h):= \gamma_1 + \gamma_2 h < 1$), then the sequence of approximate primal-dual optimizers $\{(\x_k, \blambda_k)\}_{k \in \mathbb{N}}$ generated by~\cref{alg.approx} converges linearly to an error ball around the optimal trajectory. In particular, the convergence rate is $\tau(h)$, while the asymptotical error is
\begin{subequations}\label{as_error}
\begin{align}
\limsup_{k\to\infty} \|{\blambda}_{k} \!-\! {\blambda}^*(t_{k})\| & = O\left(\frac{\varrho_{\textrm{C}}^{C}\varrho_{\textrm{P}}^{P}\,h}{1 - \tau(h)}\right) + O\left(\frac{\varrho_{\textrm{C}}^{C}\,h^2}{1 - \tau(h)}\right)\\  
\limsup_{k\to\infty} \|{\x}_{k} \!-\! {\x}^*(t_{k})\| &= O\left(\!\frac{\varrho_{\textrm{C}}^{C-1}\varrho_{\textrm{P}}^{P}\,h}{1 - \tau(h)}\!\right) \!+\! O\left(\!\frac{\varrho_{\textrm{C}}^{C-1}\,h^2}{1 - \tau(h)}\!\right)\!.
\end{align}
\end{subequations}
\end{theorem}

\begin{proof} Given in Appendix~\ref{appendix:quattro}, where the constant $\gamma_2$ is characterized as 
\begin{equation}
\gamma_2 := \frac{\kappa_f \kappa_{\A}^2}{m} \left(\frac{\kappa_f \kappa_{\A}^2 + 1}{m}C_1 C_0 + C_2\right){\varrho_{\textrm{C}}^{C-1}}(\varrho_\textrm{P}^{P}+1).
\end{equation}
And the asymptotical error is duly spelled out in terms of the problem parameters.\end{proof}

\begin{corollary}[Convergence in case of exact prediction]\label{corr.main} 
The results of Theorem~\ref{theo.main} are valid for the case of exact prediction, by letting $P \to \infty$. In particular, condition~\cref{cond:1} is verified for any $C \geq 1$. 
\end{corollary}

Theorem~\ref{theo.main} and Corollary~\ref{corr.main} dictate how the sequences generated by Algorithm~\ref{alg.approx} converge to a ball around the optimal primal-dual trajectory. For small enough sampling periods and $\tau(h)$ different enough than $1$, such that the term $1-\tau(h)$ is practically a constant for all the considered $h$, then the error ball is in the order of 
\begin{equation}\label{error_ball}
O(\varrho_{\textrm{C}}^{C}\varrho_{\textrm{P}}^{P}\, h) + O(\varrho_{\textrm{C}}^{C}\, h^2),
\end{equation}
which becomes a $O(h^2)$ error bound, every time $P$ is sufficiently large, and goes to zero if the correction step is exact ($C\to\infty$), that is every time that we solve the sampled time-invariant problems at optimality. 

The error bound $O(h^2)$, which is an improvement over a purely running scheme, for which we obtain a $O(h)$ error bound (see Corollary~\ref{corr.main}), is induced by the newly developed prediction step and it comes at the price of more restrictive conditions on $C$ and the sampling period $h$, i.e., conditions~\cref{cond:1} and \cref{cond:2}.

The parameters $P$ and $C$ need to be selected so that condition~\cref{cond:1} is satisfied. This can be achieved by computing or estimating $\varrho_{\textrm{P}}$ and $\varrho_{\textrm{C}}$ via the knowledge (or estimates) of the problem properties ($m$, $L$, $\sigma_{\max}$, $\sigma_{\min}$). Assuming that $\alpha$ and $\beta$ are chosen equal, and, e.g., $\varrho_{\textrm{P}} = \varrho_{\textrm{C}} = 0.8$, then the condition can be satisfied, e.g., with $P = 1, C\geq 5$, or $P=5, C\geq 2$, which is not as costly as it may seem. 

As can be seen from the expression of $\gamma_2$ and the condition~\cref{cond:2}, the constraint on the sampling period $h$ becomes tighter when $\gamma_2$ is large, that is when the matrix $\A$ is ill-conditioned ($\kappa_{\A}$ is large), the condition number of the problem is large ($\kappa_f$ is large), and when the time variations are important ($C_0$ and $C_2$ are large).



\section{Distributed optimization problems}\label{sec:distr}

In this section, we consider specifically distributed optimization problems. We are interested in problems of the form:
\begin{equation}\label{prob:distr}
\min_{\x\in\mathbb{R}^n} \sum_{i=1}^N f_i(\x; t),
\end{equation}
where the time-varying cost functions $f_i: \mathbb{R}^n\times \mathbb{R}_{+} \to \mathbb{R}$ verify Assumption~\ref{as.1}. In many settings, one would like to exploit the separable structure of such a cost function to decompose the optimization problem over a network of computing and communicating nodes (e.g., sensors, mobile robots). Let each node be associated with the cost function $f_i$, inducing a one-to-one mapping between nodes and local cost functions. The nodes, $i = 1, \dots, N$ can communicate via links. If two nodes $i$, $j$ share a link, we say that there is an edge connecting them. This defines a undirected graph $\mathcal{G} = (V,E)$, with vertex set $V = \{1,\dots,N\}$ and edge set $E$. The goal is now to solve~\cref{prob:distr} by allowing the nodes to communicate through their links only.  

In this case, an often employed procedure is to give each node a copy of the optimization variable, say $\y_i\in\mathbb{R}^n$ and to constrain the local variable of node $i$ to be the same as the ones of all the nodes it can communicate with. This leads to the lifted problem
\begin{equation}\label{prob:distr1}
\y^*(t):= \argmin_{\y_1\in\mathbb{R}^n,\ldots,\y_N\in\mathbb{R}^n} \sum_{i=1}^N f_i(\y_i; t),\quad \textrm{subject to}\, \A\y = {\bf 0},
\end{equation}
where $\y = (\y_1^\transp, \ldots, \y_N^\transp)^\transp \in \mathbb{R}^{nN}$ is the stacked version of all the local decision variables and $\A\in\mathbb{R}^{p \times nN}$ is a constraint matrix (i.e., the incidence matrix), whose blocks specify the fact that\footnote{We use the convention that $\y_i - \y_j = 0$, for $i\leq j$. } $\y_i = \y_j$ for all communicating couples $(i,j)$. When the underlying communication graph is connected, then the lifted problem~\cref{prob:distr} is equivalent to the original problem~\cref{prob:distr1} in the sense that each of the local optimization variable $\y_i$ at optimality is $\x^*$. 

The lifted problem~\cref{prob:distr1} is an instance of~\cref{eq.tvar}, for which the matrix $\A$ is in general rank deficient. One could reduce $\A$ to be full rank by finding a tree in the communication graph (i.e., by eliminating any linear dependent constraint), but in general one would not like to do that, since in practice convergence rates of distributed algorithms are dictated by how many links the communication graph has. The more undirected links translates in general to faster convergence.

The fact that $\A$ is rank deficient is not a problem for the proposed prediction-correction methods. However, an interesting question is whether we can perform any of the two algorithms for prediction-correction is a distributed fashion, i.e., by allowing each node $i$ to communicate only through its 1-hop communication links.

\subsection{Distributed implementation}

To obtain a distributed implementation, we require the additional assumption that:

\begin{assumption}\label{as.synch}
Communication among the nodes is synchronized; moreover, the algorithmic switching between correction and prediction is also synchronized among the nodes.  
\end{assumption}

Under Assumption~\ref{as.synch}, we claim that Algorithm~\ref{alg.approx} can be implemented on a network of communicating nodes as follows. 

\begin{claim}\label{claim:dist}
Consider the time-varying problem~\cref{prob:distr}, the communication graph $\mathcal{G} = (V,E)$, and the matrix $\A$ induced by the communication graph. Under Assumption~\ref{as.synch}, the prediction-correction Algorithm~\ref{alg.approx} can be implemented in a distributed fashion, by allowing communication only via the edge set $E$.  
\end{claim}

\begin{proof}
Given in Appendix~\ref{appendix:cinque}. 
\end{proof}

The total communication budget per time step per node (intended as the number of scalar variable transmitted) is $(P+C) N_i n$, where $N_i$ the number of neighbors of node $i$, while $P$ and $C$ are the number of prediction and correction iteration respectively.

\section{Numerical examples}\label{sec:num}

In this section, we implement our algorithm for a simple numerical example in order to assess its performance in practice. Inspired by~\cite{Xiao2006a}, we consider the following time-varying optimization problem:
\begin{equation}\label{prob-xiao}
\min_{x\in\mathbb{R}} \, \sum_{i=1}^N \underbrace{\left[\frac{1}{2}\|x- A\cos(\omega t  + \varphi_i) \|^2_2 + \log(1 + \exp(x - a_i))\right]}_{=: f_i(x; t)},
\end{equation}
with $\{a_i\}_{i=1}^N$ and $\{\varphi_i\}_{i=1}^N$ drawn from uniform probability distribution of support $[-10,10]$ and $[0,2\pi)$, respectively, while $A = 2.5$, $\omega = \pi/80$. The cost function is strongly convex and strongly smooth and $m = 1$ and $L = 1.25$, respectively. 

From a control perspective, Problem~\eqref{prob-xiao} could represent a rendezvous problem of a group of robots that would like to stay close to their moving target $A\cos(\omega t  + \varphi_i)$, and to their fixed base station located in $a_i$. Or it could represent a consensus problem, where a group of agents try to reach a compromise on their opinions on a certain matter, trading of a short-term dynamics (represented by $A\cos(\omega t  + \varphi_i)$, e.g., weekly fluctuations caused by the latest news) and a long-term one (represented by $a_i$), e.g., long-standing beliefs.   

We focus our analysis on a network of computing and communicating nodes and we fix the total number of nodes to $N = 250$, while their communication graph is randomly generated. The nodes have their local cost function $f_i(x; t)$ and they have to cooperate to determine the common decision variable $x$. By leveraging a dual decomposition approach to the described distributed optimization problem, one arrives at the problem
\begin{equation}\label{eq:numex}
\min_{y_1\in\mathbb{R}, \dots, y_N\in\mathbb{R}} \, \sum_{i=1}^N f_i(y_i; t), \quad \textrm{subject to}\, {\A}\y = {\bf 0},
\end{equation}
where $\y = (y_1, \ldots, y_N)^\transp \in \mathbb{R}^{N}$ is the stacked version of all the local decision variables and $\A$ is the constraint matrix constructed as expressed in Section~\ref{sec:distr}, and in our simulations $\kappa_{\A} = 2.39$. Problem~\cref{eq:numex} is specific version of problem~\cref{prob:distr1}, which we have analyzed in Section~\ref{sec:distr}. 

\subsection{Analysis correction-only vs. prediction-correction}

\begin{figure}
\centering
\includegraphics[width = 0.5\textwidth]{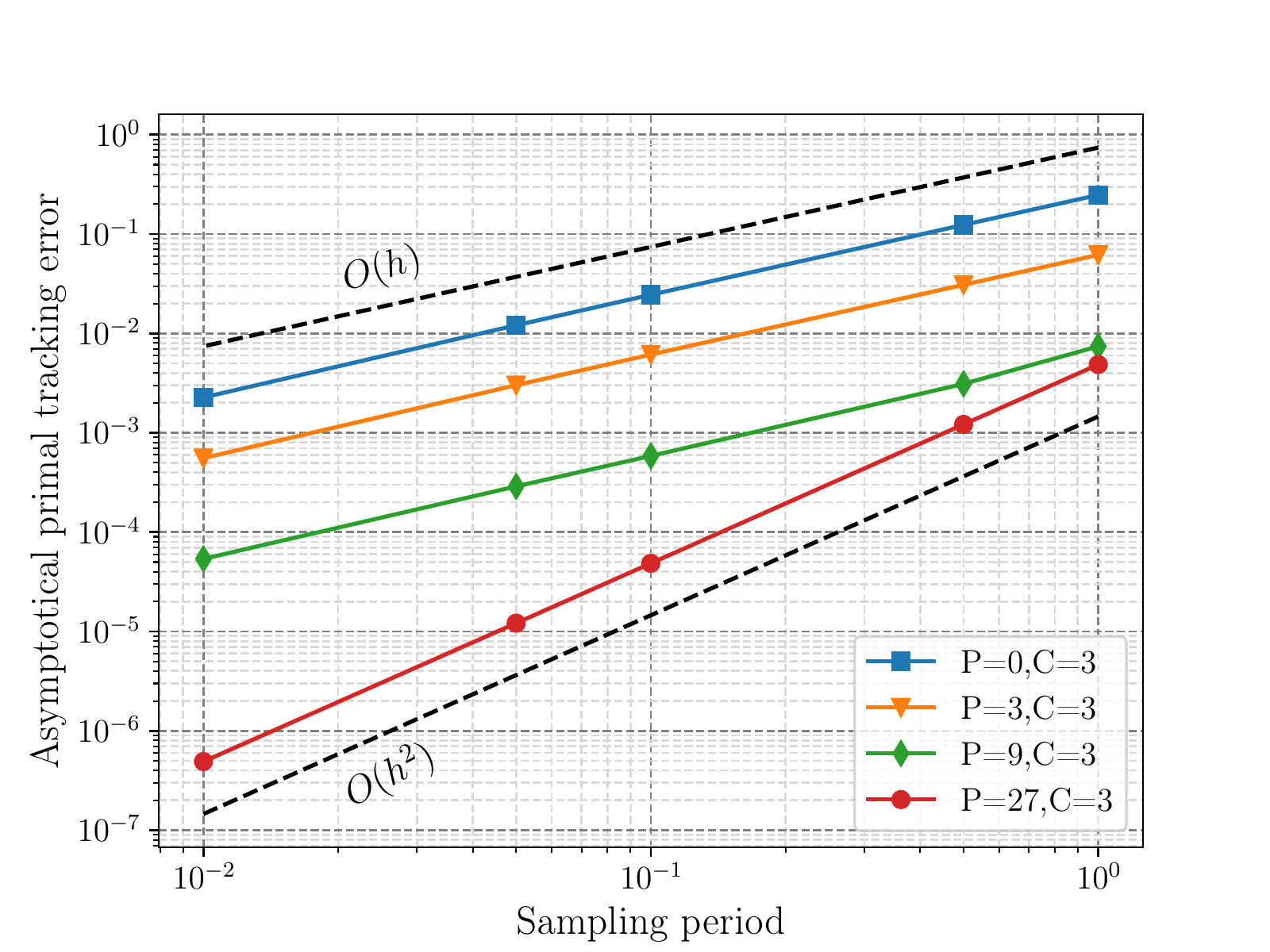}
\caption{Asymptotical tracking error as a function of sampling period for different choices number of prediction steps $P$.}
\label{fig.2sim}
\end{figure}

In the first numerical assessment, we study the proposed algorithm by varying the sampling period $h$ and for different choices of number of prediction and correction steps. In Figure~\ref{fig.2sim}, we report the results in terms of asymptotical tracking error, here computed as
\begin{equation}\label{eq:asy_err}
\max_{k\geq 5000} \{\|\y_k - \y^*(t_k)\|_2\},
\end{equation} 
whereas the final time of the simulations is $k = 10000$.

We see how a correction-only methodology (i.e., the running dual ascent discussed in Section~\ref{sec:timeinvariant}) is performing the worst, while the prediction-correction scheme with $P= 27$ (practically equivalent to the an exact prediction algorithm) is performing the best. Using a large number of prediction steps requires more computational/communication effort and therefore there is a natural trade-off between the number of prediction steps one can run and the tracking error (here captured by varying $P$). We note that, even a small number of prediction steps are beneficial in terms of asymptotical error. 

Figure~\ref{fig.2sim} depicts how the tracking error depends on the sampling period $h$ by the help of the two dashed lines, indicating a $O(h)$ and $O(h^2)$ dependence as expected by our theoretical analysis (Note that $\tau(h)$ varies less than $1\%$ for the considered sampling periods). In particular, a purely running scheme would have an asymptotical error of $O(h)$ [Cf. Corollary~\ref{co:dualascent_h}], while a prediction-correction one would have an error that approaches $O(h^2)$ when $P$ is chosen bigger and bigger [Cf. Equation~\cref{error_ball}].

We note that all the problem parameters $\{C_i\}_{i=0}^3$ can be determined and all the conditions of Theorem~\ref{theo.main} are verified. 
 
\subsection{Analysis at fixed run time}

Our second assessment regards performance of the algorithms keeping the run time per sampling period fixed, which is extremely relevant in real situations.  

Every time a new function is available, a number of correction steps are performed. The number depends on how fast we need the corrected variable to be available and the computational/communication time necessary to compute it. We fix at $r_1 h$, with $r_1<1$ the time allocated for the correction steps, while $t_{\textrm{C}}$ is the time to perform one correction step. For the above considerations, we can afford to run
\begin{equation}\label{opt.c}
C = \lfloor r_1 h/ t_{\textrm{C}}\rfloor,
\end{equation}  
correction steps. After the corrected variable is available, one can use it for the decision making process (which may require extra time to be performed). For the time-varying algorithm perspective, one can use the variable to either run $P$ gradient prediction, or $C'$ extra correction steps (to improve the corrected variable for having a better starting point when a new function becomes available). Fix at $r_2 h$, with $r_2<1$ the time allocated for the prediction (or extra correction) steps. The affordable number of prediction steps can be determined considering that $P$ prediction steps require a time equal to $\bar{t} + P t_{\textrm{P}}$, where $\bar{t}$ is the time required to evaluate once the Hessian, and time derivative of the gradient in~\cref{dual.pred.opt}, while $t_{\textrm{P}}$ is the time to perform one prediction calculation (including communication latencies). Thus, 
\begin{equation}\label{opt.p}
P = \lfloor (\,r_2 h - \bar{t}\,)/ t_{\textrm{P}}\rfloor.
\end{equation}  
The affordable extra correction steps $C'$ can be computed as in~\cref{opt.c}, substituting $r_1$ with $r_2$. 

In the simulation example, we choose $r_1 = r_2 = 0.5$, while by running the experiments on a $2.7$~GHz Intel Core i5 and by mapping the results on simple computational nodes, we empirically fix $t_{\textrm{C}} = 21$~ms, $\bar{t} = 8$~ms, $t_{\textrm{P}} = 3$~ms. Note that we have include a communication latency of $1$~ms in both $t_{\textrm{C}}$ and $t_{\textrm{P}}$, which simulates the need for communication to agree on the common decision variable, as expressed in Section~\ref{sec:distr}. Note that the correction step takes longer that a prediction step for at least two reasons. First, in the correction step one has to solve iteratively the optimization problem associated with the Lagrangian~\cref{dual.corr.opt} (here we use a Newton method), while for the prediction step such optimization problem is quadratic and unconstrained (cf.~\cref{dual.pred.opt}), so analytically solvable. Second, the aforementioned optimization problems depends on parameters (Hessian, gradient, time derivative of the gradient) that in the correction step changes for all $c \in [0,C-1]$, while they are the same for the prediction step for all $p \in [0,P-1]$ and they can be computed once.

\begin{figure}
\centering
\includegraphics[width = 0.5\textwidth]{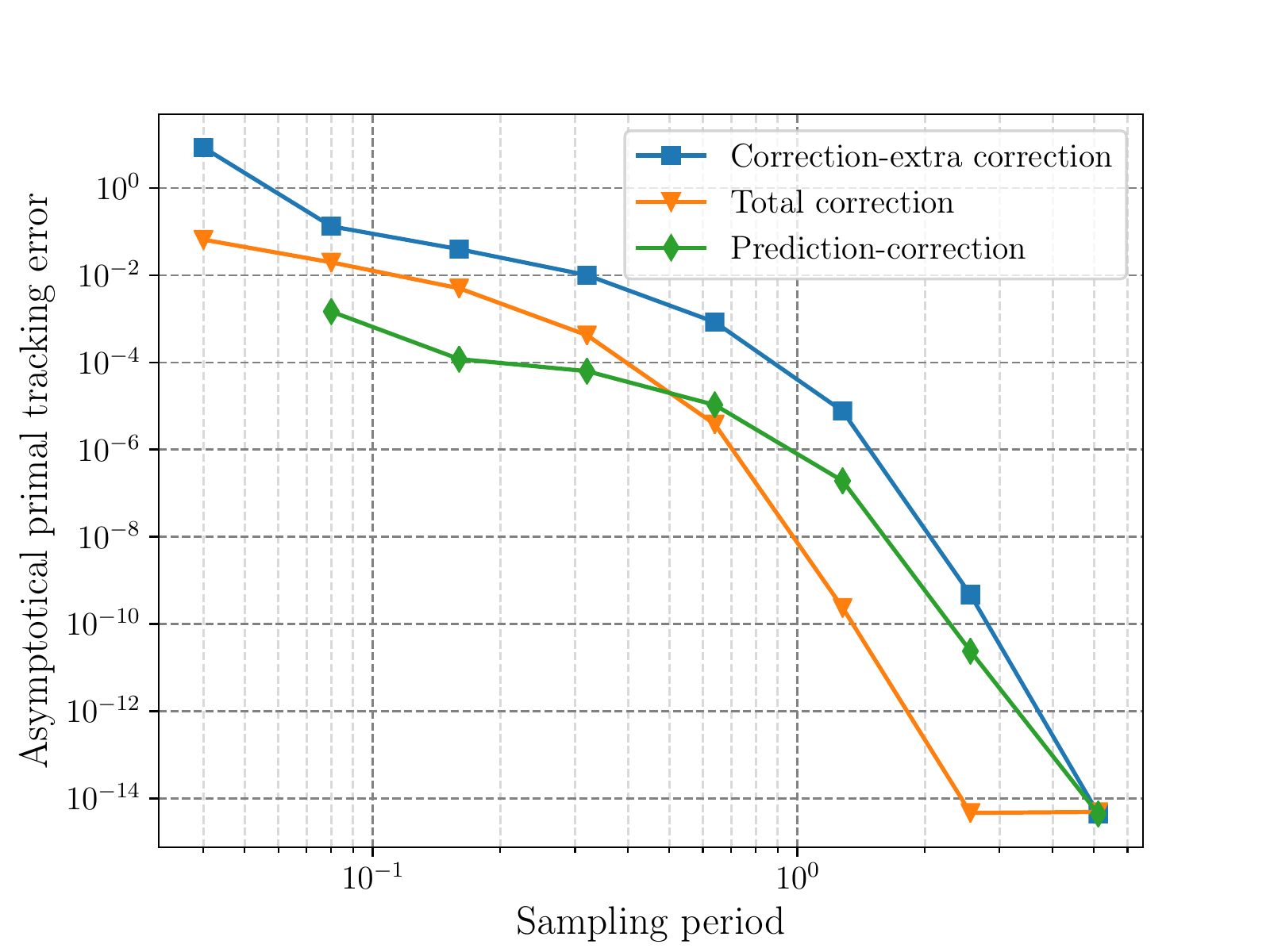}
\caption{Asymptotical tracking error as a function of sampling period for different algorithms, keeping the run time constant.}
\label{fig.3sim}
\end{figure}

In Figure~\ref{fig.3sim}, we report the results in terms of asymptotical tracking error~\cref{eq:asy_err} for the sampling period range $h \in [0.04,5.12]$~s. In the simulations, the prediction and correction steps are determined by using~\cref{opt.c} and \cref{opt.p}, so that when $h = 0.08$~s, then $P = 10$ and $C = C' = 1$, while for $h = 5.12$~s, $P = 850$ and $C = C' = 121$ (note that when $h = 0.04$~s, the prediction-correction algorithm does not satisfy the convergence assumption of Theorem~\ref{theo.main}). 

We also consider the situation in which one can use the whole sampling period to do correction, that is $r_1 = 1$, while $r_2 = 0$, and we call this case \emph{total correction}. In this case when $h = 0.08$~s, then the correction steps are $C'' = 3$, while for $h = 5.12$~s, $C'' = 243$. This total correction situation is particularly interesting when one has to make a choice whether to stop the correction steps to perform prediction, or to continue to do correction steps till a new function evaluation becomes available. Note that the correction+extra correction strategy is different from the total correction one, since the error is computed with the corrected variable (which is used for the decision making process), that is after $r_1 h$. 

The numerical results suggest that a prediction-correction strategy achieves a lower asymptotical error than performing both correction+extra correction and total correction up to a certain sampling period. This is reasonable to expect, since as $P$ and $C$ grow, the error of the prediction-correction strategy goes as $O(h^2)$, while the ones of the correction only schemes go as $O(h)$. This can be formalized as follows: the correction+extra correction strategy has an asymptotical primal error bound of
\begin{equation}\label{err.cec}
\mathrm{Err}_{\mathrm{C+EC}} = \frac{\sigma_{\max}}{m} \varrho_{\textrm{C}}^{C-1} \!\!\left(\!\!\frac{\varrho_{\textrm{C}}^{C+C'} K}{1-\varrho_{\textrm{C}}^{C}} + K\!\!\right)\! \xrightarrow[C,C'\to\infty]{}  O({\varrho_{\textrm{C}}^{C-1} h});
\end{equation}
the total correction strategy has an asymptotical primal error bound of
\begin{equation}\label{err.tc}
\mathrm{Err}_{\mathrm{TC}} = \frac{\sigma_{\max}}{m} \varrho_{\textrm{C}}^{C''-1} \left(\frac{\varrho_{\textrm{C}}^{C''} K}{1-\varrho_{\textrm{C}}^{C''}} + K\right) \xrightarrow[C''\to\infty]{}  O({\varrho_{\textrm{C}}^{C''-1} h});
\end{equation}
while the prediction-correction has an asymptotical primal error bound of
\begin{equation}\label{err.pc}
\mathrm{Err}_{\mathrm{PC}}\! =\! O\!\left(\!\frac{\varrho_{\textrm{C}}^{C-1}\varrho_{\textrm{P}}^{P}\,h}{1 - \tau(h)}\!\right) \!+\! O\left(\!\frac{\varrho_{\textrm{C}}^{C-1}\,h^2}{1 - \tau(h)}\!\right)\! \xrightarrow[P, C\to\infty]{}  O({\varrho_{\textrm{C}}^{C-1} h^2}) ;
\end{equation}
where \eqref{err.pc} is due to \eqref{as_error}, while \eqref{err.cec} and \eqref{err.tc} are generalizations of \eqref{claim2.tv} for multiple correction steps [see Appendix~\ref{appendix:sei}]. As we see, \eqref{err.cec} does not depend on $C'$ (the extra correction terms), which make these calculations superfluous, while $C''\geq 2C$ (in our case), which makes \eqref{err.tc} $<$ \eqref{err.cec}. Finally, \eqref{err.pc} is better than \eqref{err.tc} and \eqref{err.cec} for small $h$.   

The simulations indicate that, when the sampling period is small, performing prediction-correction is better than the presented alternatives, \emph{even taking into account computational and communication requirements}. In particular, \emph{(i)} w.r.t. correction+extra correction: if one has time available after the decision variable needs to be delivered and before the new cost function becomes available, doing prediction rather than extra correction appears to be the best choice; \emph{(ii)} w.r.t. total correction: it may be wise to stop the correction steps (even if one has still time before delivering the decision variable) and start the prediction ones.

\subsection{Further numerical studies}

\begin{figure*}
\centering
\includegraphics[width = 0.49\textwidth]{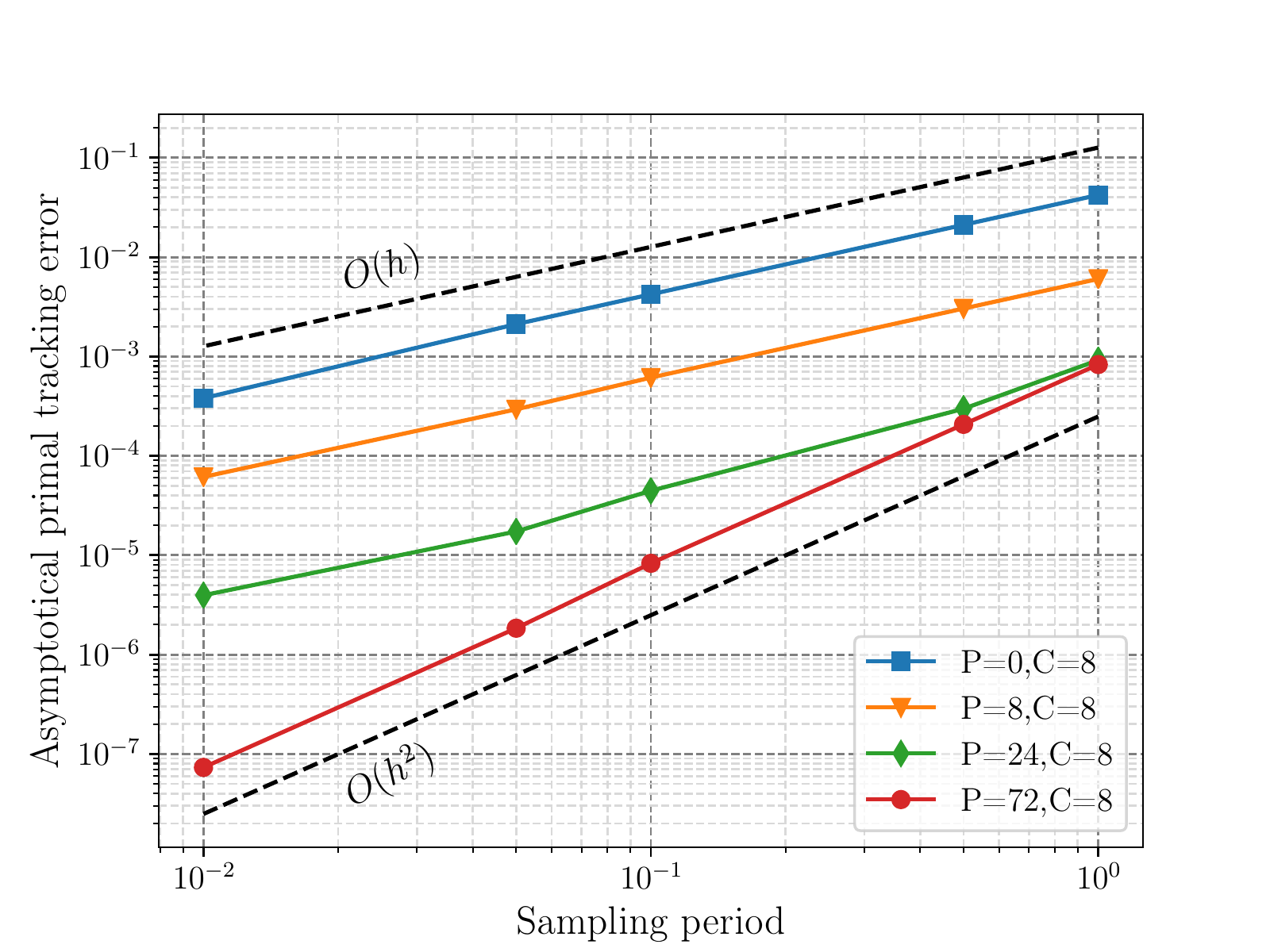}
\includegraphics[width = 0.49\textwidth]{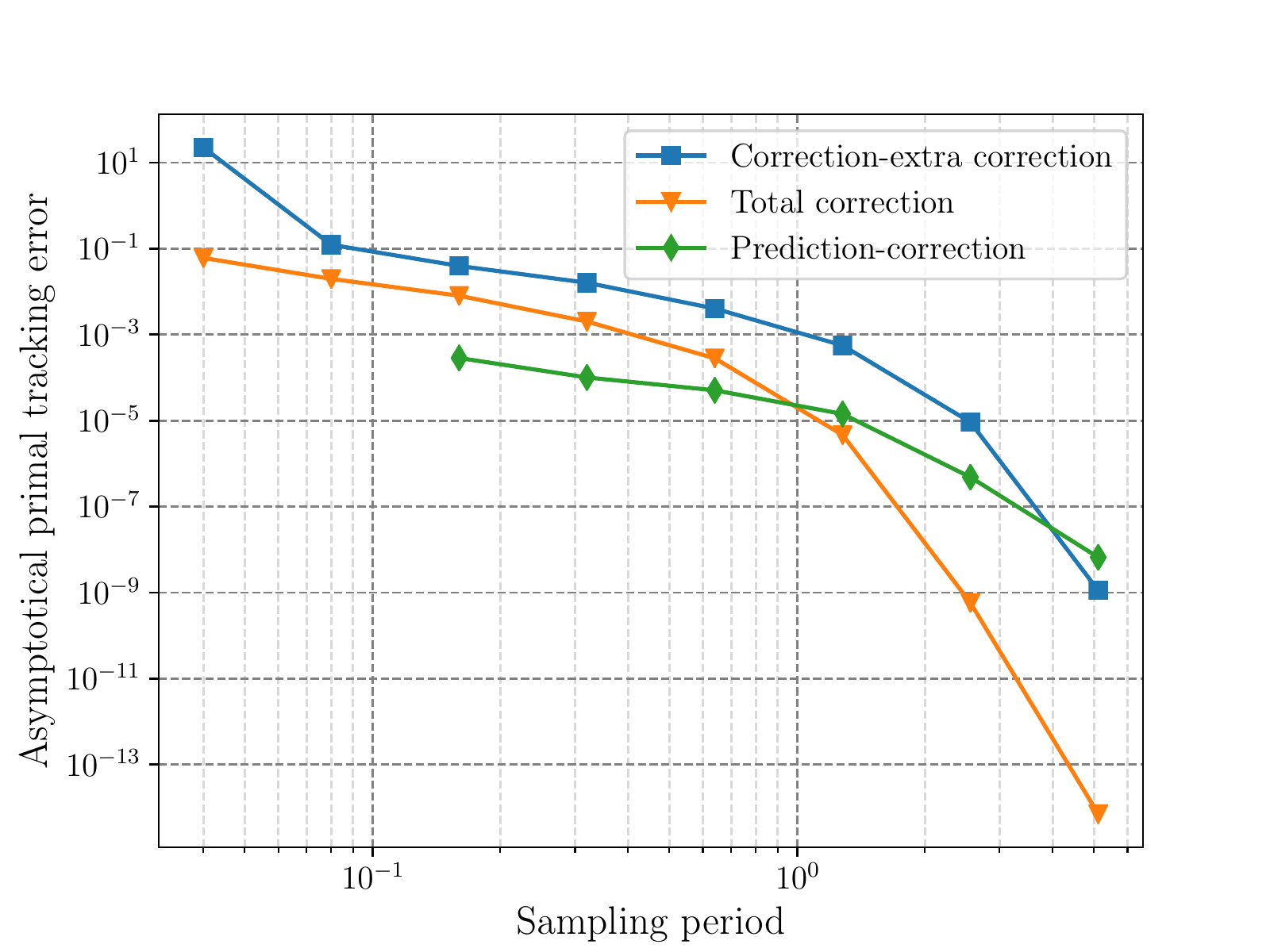}
\caption{Asymptotical tracking error performance for $k_f = 3.25$ and other parameters left the same.}
\label{fig.4sim}
\end{figure*}

\begin{figure*}
\centering
\includegraphics[width = 0.49\textwidth]{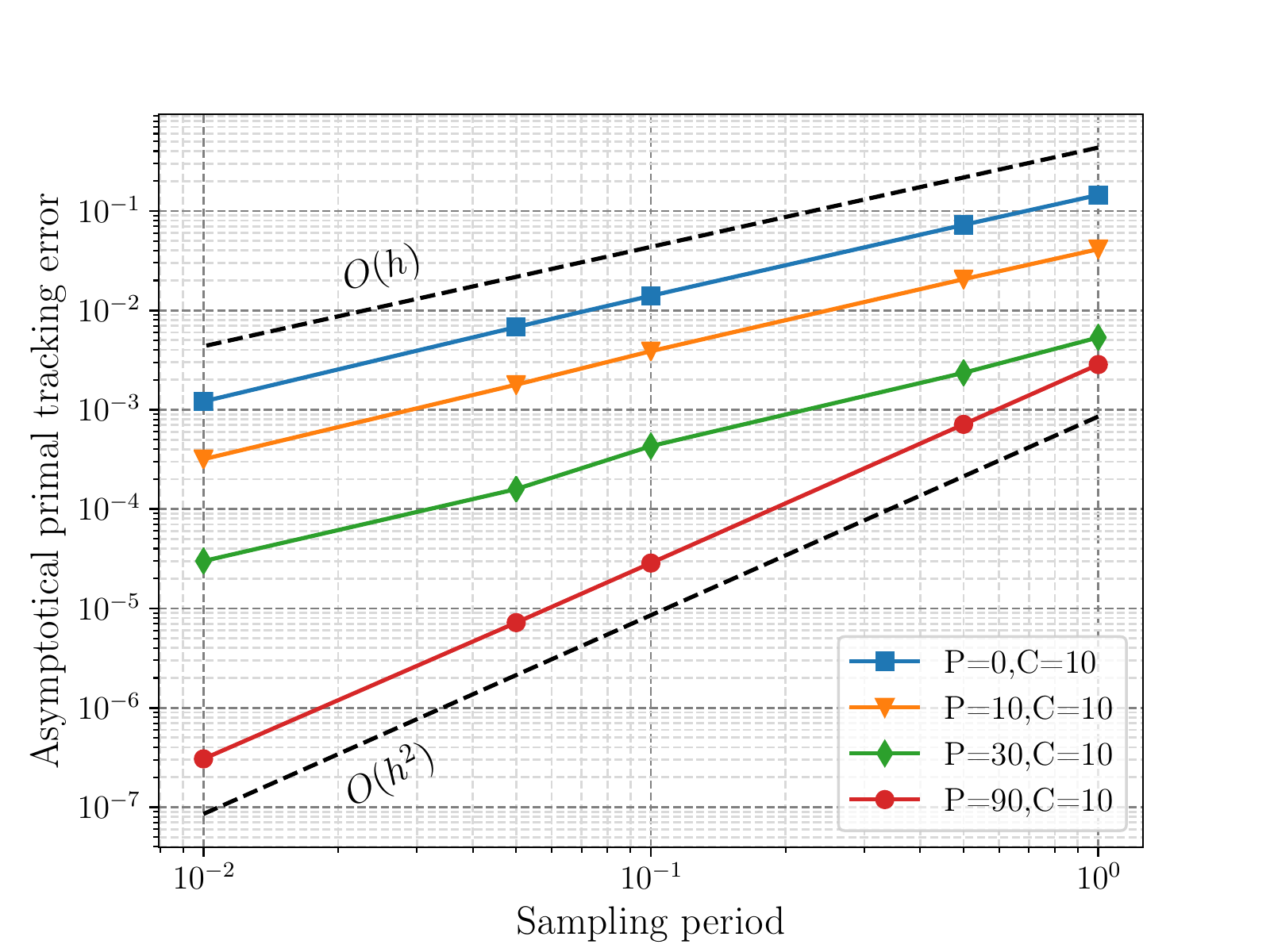}
\includegraphics[width = 0.49\textwidth]{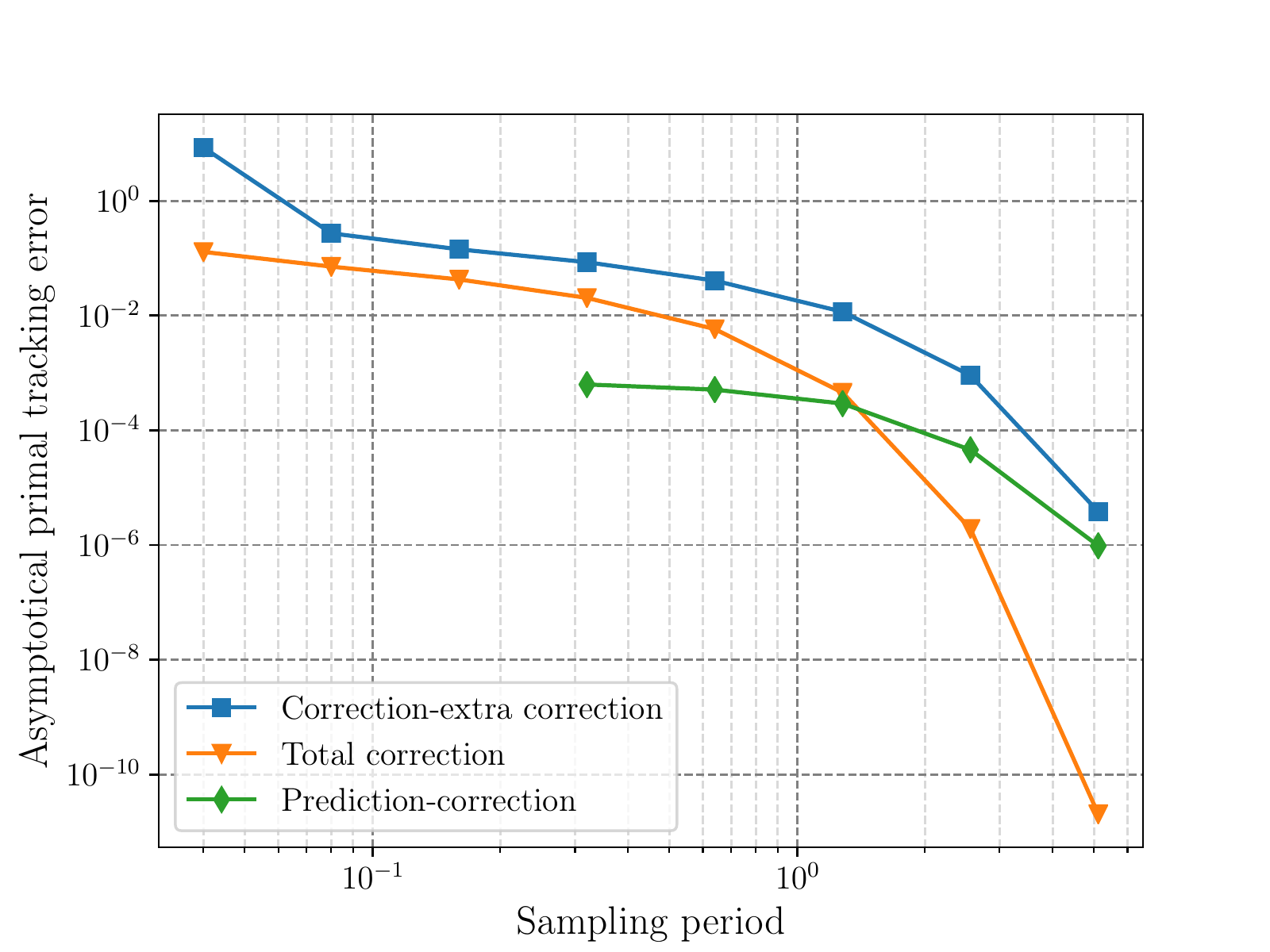}
\caption{Asymptotical tracking error performance for $k_{\A} = 4.28$ and other parameters left the same.}
\label{fig.5sim}
\end{figure*}

\begin{figure*}
\centering
\includegraphics[width = 0.49\textwidth]{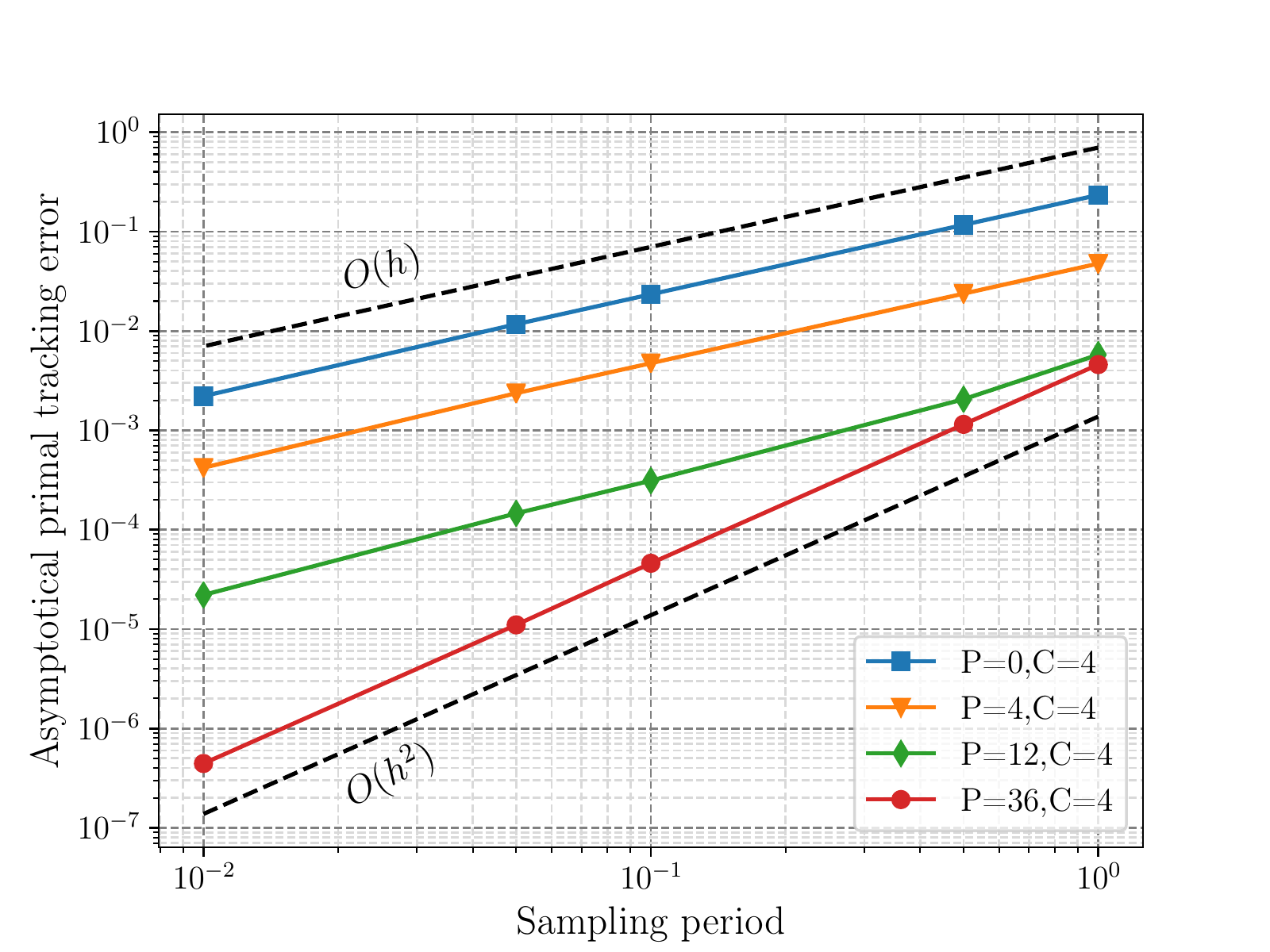}
\includegraphics[width = 0.49\textwidth]{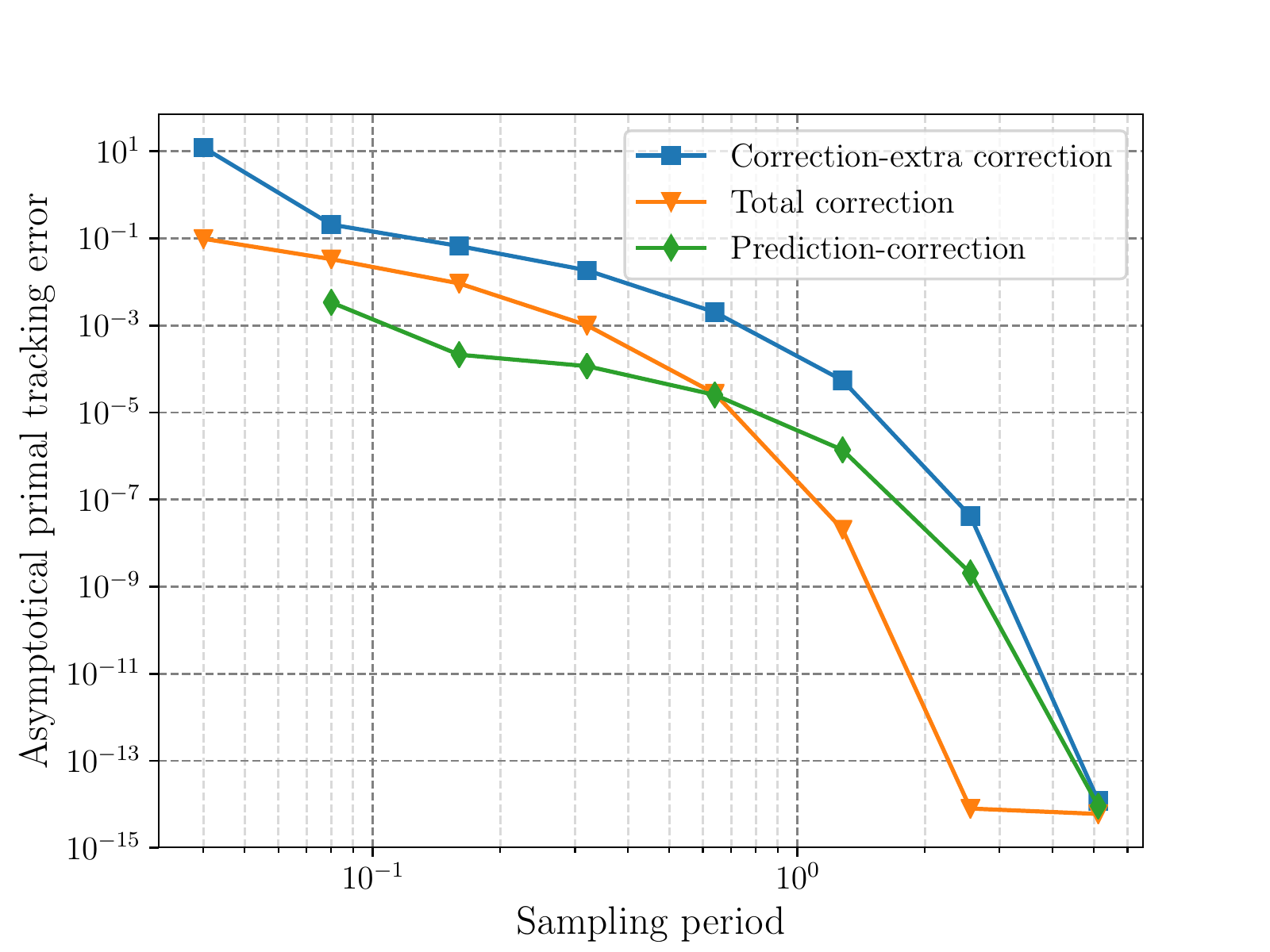}
\caption{Asymptotical tracking error performance for $N = 500$, $\kappa_{\A} = 2.54$ and other parameters left the same.}
\label{fig.6sim}
\end{figure*}

We report here further numerical studies which are qualitatively very similar to the ones just presented. In particular, we report that both \emph{(i)} changing the condition number of the function $\kappa_f$ from $1.25$ to $3.25$ [Figure~\ref{fig.4sim}] and \emph{(ii)} changing the condition number of the incidence matrix $\kappa_{\A}$ from $2.39$ to $4.28$ [Figure~\ref{fig.5sim}], require more prediction and correction steps to achieve the same asymptotical error bounds; whereas \emph{(iii)} increasing the number of nodes from $N=250$ to $N=500$ (while having $\kappa_{\A} = 2.54$) [Figure~\ref{fig.6sim}], has very limited effect in the number of prediction and correction steps required. 

\section{Conclusions}\label{sec:conc}

We have developed dual prediction-correction methods to track the solution trajectory of time-varying linearly constrained convex programs. The proposed methods have a better theoretical and numerical performance with respect to more classical strategies. We have characterized the convergence properties and asymptotical tracking error of all the methods and shown how the error depends on the problem instance parameters and sampling period. 



\newpage
\appendices

\section{Proof of Proposition~\ref{prop:properties}}\label{appendix:uno}

\begin{proof}
Call $f_k : = f(\cdot;t_k)$. The primal optimizer of~\cref{eq.tinv} $\x^*(t_k)$ is unique since $f_k$ is strongly convex. Examine the optimality condition,
\begin{equation}\label{eq.opt1}
\nabla_{\x} f_k(\x^*(t_k)) + \A^\transp \blambda^*(t_k) = {\bf 0}.
\end{equation}
By strong smoothness $\nabla_{\x} f_k(\x^*(t_k))$ is unique. In fact, if there were two distinct $\nabla_{\x} f_k(\x_1)$ and $\nabla_{\x} f_k(\x_2)$, for the same $\x_1 = \x_2$, then one could derive a contradiction by using the strong smooth inequality
\begin{equation}
\|\nabla_{\x} f_k(\x_1) - \nabla_{\x} f_k(\x_2) \|\leq L \|\x_1 - \x_2\|.
\end{equation}
Thus, $\A^\transp \blambda^*(t_k)$ is unique.

By the fundamental theorem of linear algebra (or alternatively, Fredholm alternative theorem) \cite{Strang2016}, $\blambda^*(t_k)$ can be decomposed in two parts as $\blambda^*(t_k) = \blambda^*_0(t_k) + \v_k$, for which $\blambda^*_0(t_k) \in \im(\A)$ and $\v_k \in \nul(\A^\transp)$. 

In the full rank case, the nullspace of $\A^\transp$ is void and $\blambda^*(t_k) \in \im(\A)$. In the rank deficient case, we only concentrate on $\blambda^*_0(t_k) \in \im(\A)$. Uniqueness of $\blambda^*_0(t_k)$ is proven by contradiction: assume that $\blambda^*_0(t_k)$ is not unique and one has two variables $\blambda_1 \neq \blambda_2$ for which $\A^\transp \blambda_1  = \A^\transp \blambda_2$. Since both variables lie in the image of $\A$ (and not in the nullspace of $\A^\transp$), it has to be $\A^\transp \blambda_1 \neq 0$, $\A^\transp \blambda_2 \neq 0$ as well as for any of their linear combinations. Therefore, it has to be $\A^\transp (\blambda_1 - \blambda_2) \neq 0$ for all $\blambda_1 \neq \blambda_2$, from which a contradiction arises. Therefore $\blambda^*_0(t_k)$ must be unique.
\end{proof}

\section{Proof of Theorem~\ref{theo.dualascent}}\label{appendix:due}

\begin{proof}
The proof is reported here for completeness, it can be found e.g. in~\cite{Ryu2015,Simonetto20XX}. Call $f_k : = f(\cdot;t_k)$. The proof relies on the properties of the conjugate function of $f_k$, defined as $f^{\star}_k(\y) := \sup_{\x\in\mathbb{R}^n}\{\y^\transp \x - f_k(\x)\}$, and on the properties of the differential operator $\partial g$ of a convex function $g$. In particular, if $f_k$ is strongly convex for all $\x\in\mathbb{R}^n$ with parameter $m$, then $f^{\star}_k$ is strongly smooth with parameter $1/m$ for all $\y \in \mathbb{R}^n$, while if $f_k$ is strongly smooth with parameter $L$ for all $\x\in\mathbb{R}^n$, then $f^{\star}_k$ is strongly convex with parameter $1/L$ for all $\y \in \mathbb{R}^n$. Furthermore, for the differential operators of $f_k$ and $f^{\star}_k$, one has $\partial f_k^{-1} = \partial f^{\star}_k$. 

Consider now~\cref{dual.ascent1}, which can be written in terms of optimality condition as 
\begin{equation}
\partial f_k(\x_{i+1}) + \A^\transp \blambda_i = {\bf 0} \iff  \x_{i+1} = \partial f^{\star}_k(-\A^\transp \blambda_i),
\end{equation}
where we have used the identity $\partial f_k = \nabla f_k$, since $f_k$ is differentiable. The dual function $q_k(\blambda) : = \min_{\x \in \mathbb{R}^n}\{f_k(\x) + \blambda_{i}^\transp (\A\x-\b) \}$ has gradient,
\begin{equation}\label{dual:fun}
\partial q_k(\blambda) = \A\x_{i+1} - \b = \A \partial f^{\star}_k(-\A^\transp \blambda_i) -\b.
\end{equation}

{\bf Full row rank $\A$.} Due to~\cref{dual:fun}, the negative of the dual function $-q_k$ is strongly smooth with constant $\sigma_{\max}^2/m$ and strongly convex with constant $\sigma_{\min}^2/L$. The dual ascent~\cref{dual.ascent} is a dual gradient iteration on $-q_k$ and it converges for all $\alpha < 2 m/\sigma_{\max}^2$, with linear convergence rate $\varrho = \max\{|1-\alpha L/\sigma_{\min}^2|,|1-\alpha m/\sigma_{\max}^2|\}$. Therefore, 
\begin{equation}
\|\blambda_{i+1} - \blambda^*(t_k)\| \leq \varrho \|\blambda_{i} - \blambda^*(t_k)\|, \end{equation}
which is claim~\cref{claim1}. By using the optimality condition for~\cref{dual.ascent1}, 
\begin{equation}
\nabla_{\x} f_k(\x_{i+1}) + \A^\transp \blambda_i = \nabla_{\x} f_k(\x^*(t_k)) + \A^\transp \blambda^*(t_k) = {\bf 0}.
\end{equation}
By algebraic manipulations and by using strong convexity,
\begin{multline}
m \|\x_{i+1}- \x^*(t_k)\| \leq \|\nabla_{\x} f_k(\x_{i+1})-\nabla_{\x} f_k(\x^*(t_k))\| \\ =\| \A^\transp (\blambda_i- \blambda^*(t_k))\| \leq \sigma_{\max} \|\blambda_i- \blambda^*(t_k)\|,
\end{multline}
from which claim~\cref{claim2}. 

{\bf Rank deficient $\A$.} To prove the contraction property in this case, we only need to re-work the strong convexity property of $-q_k$, since now $\sigma_1 = 0$.  To do that, we need to show that the functions $-q_k$ have a strong convex-like property for all $\blambda \in \im(\A)$ and that the iterations~\cref{dual.ascent} generates $\blambda_i \in \im(\A)$ (i.e., keeps the dual variable feasible). The second claim is easy to show since $\blambda_0 \in \im(\A), \b \in \im(\A)$ and 
\begin{equation}
\blambda_{i+1} = \blambda_{i} + \alpha (\A \x_{i+1}-\b) \in \im(\A).  
\end{equation}
To show the first claim, we recall that $\partial q_k(\blambda) = \A\partial \cf_k (-\A^\transp \blambda) - \b$ (as proved in the proof of full rank $\A$). Therefore, for all $\blambda,\bmu\in\im(\A)$:
\begin{multline}\label{eq.dummy00}
(\partial q_k(\blambda) - \partial q_k(\bmu))^\transp (\bmu - \blambda) = \\ =(\partial \cf_{k} (-\A^\transp \blambda) - \partial \cf_{k} (-\A^\transp \bmu))^\transp \A^\transp (\bmu - \blambda) \geq \\ \geq  \sigma_{\min}^2/L \|\blambda - \bmu\|^2,  
\end{multline}
where the last inequality comes from the fact that $\blambda,\bmu\in\im(\A)$ and by the fact that $\A^\transp (\bmu - \blambda) = {\bf 0}$ iff $\bmu = \blambda$, for the fundamental theorem of linear algebra~\cite{Strang2016}. Result~\cref{eq.dummy00} implies strong monotonicity of $-\partial q_k(\blambda)$ for all $\blambda \in \im (\A)$, and therefore strong convexity of $-q_k(\blambda)$ for all $\blambda \in \im (\A)$.  Then the contraction property follows from the fact that $-q_k$ is both strongly smooth with constant $\sigma^2_{\max}/m_k$ as easy to show, and strongly convex (over the restricted domain). The rest follows as in the proof of the full row rank case. \end{proof}

\section{Proof of Theorem~\ref{theo.implicit}}\label{appendix:tre}

In order to prove Theorem~\ref{theo.implicit}, we need a general result on quadratic programs of a special form. 

\begin{proposition}\label{prop.perturbation}
Consider the strongly convex quadratic program, 
\begin{equation}\label{eq.qpert}
\min_{\delta\x \in \mathbb{R}^n} \frac{1}{2} \delta\x^\transp \Q \delta\x + \c^\transp \delta\x, \quad \textrm{subject to } \A \delta \x = {\bf 0}, 
\end{equation}
with unique primal-dual optimizers $(\delta\x^*, \delta\blambda^*\in\im(\A))$. Let the eigenvalues of $\Q$ be lower bounded by $m$ and upper bounded by $L$, while let the matrix $\A$ have the singular values ordered as in Section~\ref{sec:cases}. Then,
\begin{equation}\label{eq.perturb}
\|\delta \x^*\| \leq \frac{1}{m}\left(1 + \frac{L}{m}\frac{\sigma_{\max}^2}{\sigma_{\min}^2}\right) \|\c\|, \qquad \|\delta \blambda^*\| \leq \frac{L}{m}\frac{\sigma_{\max}}{\sigma_{\min}^2} \|\c\|.
\end{equation}
\end{proposition}

\begin{proof}
The optimality condition for~\cref{eq.qpert}: $\Q\delta\x^* +\c + \A^\transp \delta \blambda^* = {\bf 0}$ yields,
\begin{equation}\label{dummy.opt0}
\delta\x^* = - \Q^{-1}(\c + \A^\transp \delta \blambda^*).
\end{equation}
The dual problem of~\cref{eq.qpert} reads, 
\begin{equation}\label{dualqp1}
\min_{\delta\blambda}\quad  \frac{1}{2}\delta\blambda^\transp \A\Q^{-1}\A^\transp \delta\blambda +  \c^\transp \Q^{-1} \A^\transp \delta\blambda,
\end{equation}
whose optimality condition reads,
\begin{equation}\label{dummy.opt}
\A\Q^{-1}\A^\transp \delta\blambda^* = - \A \Q^{-1} \c. 
\end{equation}
If $\A$ is full row rank, then $\delta\blambda^*$ is unique and $\delta\blambda^*\in\im(\A)$, otherwise there exists a unique $\delta\blambda^*\in\im(\A)$ (see Proposition~\ref{prop:properties}). We focus on the unique $\delta\blambda^*\in\im(\A)$. In this case $\A^\transp \delta\blambda^* \neq 0$ if $\delta \blambda^* \neq 0$ and therefore we can multiply both sides of \cref{dummy.opt} by $\delta \blambda^{*,\transp}$, obtaining,
\begin{equation}\label{dummy.opt1}
\delta\blambda^{*,\transp} \A\Q^{-1}\A^\transp \delta\blambda^* = - \delta\blambda^{*,\transp}\A \Q^{-1} \c. 
\end{equation}
Bounding,
\begin{multline}\label{dummy.opt1_1}
\frac{\sigma_{\min}^2}{L} \|\delta\blambda^*\|^2 \leq \|\delta\blambda^{*,\transp} \A\Q^{-1}\A^\transp \delta\blambda^*\| =\\ \| \delta\blambda^{*,\transp}\A \Q^{-1} \c\| \leq \| \delta\blambda^*\| \|\A\| \| \Q^{-1}\| \| \c\| \leq  \|\delta\blambda^*\| \frac{\sigma_{\max}}{m} \|\c\|, 
\end{multline}
and finally, by dividing by the non-zero, finite $\|\delta \blambda^*\|$, one derives the claim~\cref{eq.perturb} on the dual variable. For the claim~\cref{eq.perturb} on the primal variable, one uses~\cref{dummy.opt0}, which can be upper bounded as
\begin{equation}\label{xx}
\|\delta \x^*\| \leq \|\Q^{-1}\|(\|\c\| + \|\A^\transp\|\|\delta\blambda^*\|) \leq \frac{1}{m}\left(1 + \frac{L}{m}\frac{\sigma_{\max}^2}{\sigma_{\min}^2} \right)\|\c\|.
\end{equation}
\end{proof}

We are now ready for the proof of Theorem~\ref{theo.implicit}.

\begin{proof}
To determine the bounds in~\cref{bounds.xl}, we use a Taylor expansion. In particular, call $\Q = \nabla_{\x\x} f(\x^*(t_k);t_k)$ and $\c =h \nabla_{t\x} f(\x^*(t_k);t_k)$. Then, $\x^*(t_{k+1})$ can be computed as the solution of
\begin{multline}\label{hot}
\nabla_{\x}f(\x; t_{k+1}) + \A^\transp \blambda = \Q (\x - \x^*(t_k)) + \c + \textrm{h.o.t.} +\\ + \A^\transp (\blambda-\blambda^*(t_k)) = {\bf 0},
\end{multline}
where h.o.t. stands for the higher order terms of the expansion. The results provided in~\cref{bounds.xl} will be valid when the higher order terms are negligible with respect to the leading terms (i.e., locally), or when $C_1 = C_2 = C_3 = 0$, i.e., when the higher order terms are identically zero. Problem~\cref{hot} can be put in the form of~\eqref{eq.qpert} by neglecting the h.o.t., and in particular, its solution is $\|\delta \x^*\| = \|\x^*(t_{k+1}) - \x^*(t_{k})\|$ and $\|\delta \blambda^*\| = \|\blambda^*(t_{k+1}) - \blambda^*(t_{k})\|$. By the bounds on $\Q = \nabla_{\x\x} f(\x^*(t_k);t_k)$, the upper bound on $\|\c\| = \|h \nabla_{t\x} f(\x^*(t_k);t_k)\| \leq C_0 \, h$, and by using Proposition~\ref{prop.perturbation}, ones derives the claims~\cref{bounds.xl.primal} and~\cref{bounds.xl.dual}. 
\end{proof}








\section{Proof of Theorem~\ref{theo.main}} \label{appendix:quattro}

\subsection{Preliminaries and definitions}

We begin the convergence analysis by deriving an upper bound on the norm of the approximation error incurred by the Taylor expansion in~\cref{taylor}. In particular, given  the optimal primal-dual solutions $(\x^*(t_k), \blambda^*(t_k)\in\im(\A))$ and $(\x^*(t_{k+1}), \blambda^*(t_{k+1})\in\im(\A))$ at $t_k$ and $t_{k+1}$, respectively, compute the optimal prediction step via the Taylor approximation~\cref{taylor} and indicate the optimal prediction as $(\x^*_{k+1|k}= \x^*(t_k)+\delta\x, \blambda^*_{k+1|k}= \blambda^*(t_k)+\delta\blambda)$. The objective is to bound the error:
\begin{align}\label{eq.def.delta}
{e}_k : = \max\{\|\x^*_{k+1|k} - \x^*(t_{k+1})\|, \|\blambda^*_{k+1|k} - \blambda^*(t_{k+1})\|\},
\end{align}
which is committed when the optimal couple $(\x^*(t_{k+1}), \blambda^*(t_{k+1}))$ is replaced by the predicted one $(\x^*_{k+1|k},\blambda^*_{k+1|k})$.  

To ease notation, we define the following problem specific quantities:
\begin{subequations}\label{definition}
\begin{align}
\Delta_1 &:= \frac{\kappa_f \kappa_{\A}^2 + 1}{m}, \qquad \Delta_2 := \frac{\kappa_f \kappa_{\A}}{\sigma_{\min}},\\
\Delta_3 &:= C_1 C_0^2 \Delta_1^2 /2 +   \Delta_1 C_2 C_0 + C_3/2, \\ \Delta_4 &:= \Delta_1 C_1 C_0 + C_2.
\end{align}
\end{subequations}

\vspace{1mm}
\begin{proposition}\label{prop.err}
Let Assumptions~\ref{as.1}, \ref{as.3}, and \ref{as.4} hold. The following holds:
\begin{subequations}
\begin{align}
\|\x^*_{k+1|k} - \x^*(t_{k+1})\|&\leq \Delta_1 \Delta_3\, h^2 \\
\|\blambda^*_{k+1|k} - \blambda^*(t_{k+1})\| &\leq \Delta_2 \Delta_3\, h^2,
\end{align}
\end{subequations}
where $\Delta_1$, $\Delta_2$, and $\Delta_3$ are defined in~\cref{definition}. Thus, the error $e_k$ is upper bounded as
\begin{equation}
e_k \leq \max\{\Delta_1, \Delta_2\}\, \Delta_3 \, h^2 = O(h^2). 
\end{equation}
\end{proposition}

\begin{proof}
Let us start by simplifying the notation. Define
\begin{subequations}\label{simplif}
\begin{align}
\nabla_{\x} f_i \!= \!\nabla_{\x} f (\x^*(t_{k+i}); t_{k+i}), \quad & \Q_i \!= \!\nabla_{\x\x} f (\x^*(t_{k+i}); t_{k+i}) \\ \label{simplif22}
\c_i = \nabla_{t\x} f (\x^*(t_{k+i}); t_{k+i}),\quad & \x_i = \x^*(t_{k+i}), \, \x = \x^*_{k+1|k},\\
& \blambda_i = \blambda^*(t_{k+i}),\, \blambda = \blambda^*_{k+1|k}.
\end{align}
\end{subequations}
With this notation in place, ${e}_k =\max\{\|\x - \x_{1}\|,\|\blambda - \blambda_{1}\|\}$ [Cf.~\cref{eq.def.delta}]. In addition, $(\x,\blambda)$ is computed by the optimal conditions [Cf.~\cref{taylor}]
\begin{equation}\label{x.l}
\nabla_{\x} f_0 + \Q_0 (\x - \x_0) + h\, \c_0 + \A^\transp \blambda = \mathbf{0}, \qquad
\A \x = \b.
\end{equation}
while $\x_1$ is the solution of 
\begin{equation}
\nabla_{\x} f_1 + \A^\transp \blambda_1 = \mathbf{0}, \qquad
\A \x_1 = \b.
\end{equation}


Consider the solution mapping:
\begin{multline}
s(\p) :=\Big\{\y,\bmu\in\im(\A)\Big| \\ \begin{array}{c}\nabla_{\x} f_0 + \Q_0 (\y- \x_0) + h\, \c_0 +\A^\transp\bmu + \p = \mathbf{0}\\ \A \y=\b \end{array} \Big\}.
\end{multline}
The mapping $s(\p)$ is every-where single-valued (due to Proposition~\ref{prop:properties}), while for any two values of the parameter $\p$, say $\p'$ and $\p''$, then, 
\begin{subequations}\label{p.pert}
\begin{align}
\hspace*{-.25cm} \Q_0 (\y(\p')\!- \!\y(\p'')) \!+\! \A^\transp(\bmu(\p')\!-\!\bmu(\p'')) \!+\! \p'\!-\!\p'' & \!=\! \mathbf{0}, \\ \A (\y(\p') - \y(\p'')) &\! =\!\mathbf{0}.
\end{align}
\end{subequations}
By using Proposition~\ref{prop.perturbation} on \cref{p.pert} with $\c = \p'-\p''$, 
\begin{subequations}\label{bounds.p}
\begin{align}
\|\y(\p') - \y(\p'')\| &\leq \left(\frac{L}{m}\frac{\sigma_{\max}^2}{\sigma_{\min}^2} + 1 \right)\frac{1}{m}\, \|\p'-\p''\| =  \nonumber \\ & \hspace*{2.5cm}= \Delta_1  \|\p'-\p''\|, \\
\|\bmu(\p') - \bmu(\p'')\| &\leq  \frac{L}{m}\frac{\sigma_{\max}}{\sigma_{\min}^2} \|\p' - \p''\|= \Delta_2  \|\p'-\p''\|.
\end{align}
\end{subequations}

Let $\p' = {\bf 0}$ and $\p'' = \nabla_{\x} f_1-(\nabla_{\x} f_0 + \Q_0 (\x_1- \x_0) + h\, \c_0)$, one obtains $\y(\p') = \x$, $\bmu(\p') = \blambda$, and $\y(\p'')=\x_1$, and $\bmu(\p'')=\blambda_1$, which means, 
\begin{multline}\label{eq-dummy}
e_k = \max\{\|\x - \x_{1}\|,\|\blambda - \blambda_{1}\|\} \leq  \max\left\{\Delta_1, \Delta_2\right\} \times \\ \times \|\nabla_{\x} f_0 - \nabla_{\x} f_1  + \Q_0 (\x_1- \x_0) + h\, \c_0\|.
\end{multline}
Consider now the right-hand-side of~\cref{eq-dummy}: it is nothing else but the error of the truncated Taylor expansion of $\nabla_{\x} f_1$:
\begin{equation}
\nabla_{\x} f_1 - \nabla_{\x} f_0 =  \Q_0 (\x_1-\x_0) + h \c_0 +\mathbold{\epsilon},
\end{equation} 
where the error $\mathbold{\epsilon}$ can be bounded as
\begin{multline}
\|\mathbold{\epsilon}\| \leq \frac{1}{2}\Big(\|\nabla_{\x\x\x} f\| \|\x_1- \x_0\|^2 + h\, \|\nabla_{t\x\x} f\| \|\x_1- \x_0\| +\\  h\, \|\nabla_{\x t\x} f\| \|\x_1- \x_0\|  + {h^2} \|\nabla_{t t\x} f\|\Big),
\end{multline}
and by using the upper bounds in Assumption~\ref{as.4}, 
\begin{multline}\label{perdopo1}
\|\nabla_{\x} f_0 - \nabla_{\x} f_1  + \Q_0 (\x_1- \x_0) + h\, \c_0\| \leq \\ \leq \frac{1}{2}\, C_1 \|\x_1 - \x_0\|^2 + h\, C_2 \|\x_1 - \x_0\| + \frac{1}{2} h^2\, C_3. 
\end{multline}
By using the bound~\cref{bounds.xl} on the variability of the optimizers $\x_1$ and $\x_0$, then 
\begin{equation}
\|\p'-\p''\|\!\leq\! \left(\frac{1}{2} h^2\,C_1 {\Delta_1^2 C_0^2} \!+ \!h^2\, C_2 {\Delta_1 C_0} \!+\! \frac{1}{2} h^2\, C_3\right) \!= \!\Delta_3 \, h^2, 
\end{equation}
which, by substituting into~\cref{bounds.p}, proves Proposition~\ref{prop.err}. 
\end{proof}

We then look at the optimal prediction error, i.e., the distance between the exact predicted pair $(\x_{k+1|k},\blambda_{k+1|k})$ and the primal-dual optimizer at time step $t_{k+1}$, $(\x^*(t_{k+1}),$ $\blambda^*(t_{k+1})\in\im(\A))$ can be bounded as the following proposition. 

\begin{proposition}\label{prop.intheo.1}
Under the same assumptions and notation of Theorem~\ref{theo.main}, let $(\x_{k+1|k}, \blambda_{k+1|k})$ be the exact predicted step obtaining by solving~\cref{qp1} at optimality. Let $\Delta_1, \Delta_2, \Delta_3, \Delta_4$ be defined as in~\cref{definition}. We have that
\begin{subequations}\label{eq.inside}
\begin{align}
\|\x_{k+1|k} - \x^*(t_{k+1})\| & \leq (1 + h\, \Delta_1 \Delta_4) \|\x_k-\x^*(t_k)\| + \nonumber \\ & \hspace*{2cm}+ \Delta_1 \Delta_3\, h^2,\\
\|\blambda_{k+1|k} - \blambda^*(t_{k+1})\| & \leq \|\blambda_k - \blambda^*(t_{k})\| + \nonumber \\ & \hspace*{-1cm}+ h\, \Delta_2 \Delta_4 \|\x_k-\x^*(t_k)\|+ \Delta_2 \Delta_3\, h^2.
\end{align}
\end{subequations}
\end{proposition}

\begin{proof}
We proceed as in the proof of Proposition~\ref{prop.err}. We use similar simplifications of~\cref{simplif}, as
\begin{subequations}\label{simplif2}
\begin{align}
\nabla_{\x} f_k = \nabla_{\x} f (\x_k; t_{k}), \quad & \Q_k = \nabla_{\x\x} f (\x_k; t_{k}) \\
\c_k = \nabla_{t\x} f (\x_k; t_{k}),\quad & \x = \x_{k+1|k}, \quad \blambda =  \blambda_{k+1|k}.
\end{align}
\end{subequations}
while $\nabla_{\x} f_1$, $\x_1$, and $\blambda_1$, $\Q_0$, $\c_0$, $\nabla_{\x} f_0$, $\x_0$, and $\blambda_0$ are defined just as in~\cref{simplif}. The error $\|\x_{k+1|k} - \x^*(t_{k+1})\|$ is now $\|\x - \x_1\|$, while $\|\blambda_{k+1|k} - \blambda^*(t_{k+1})\|$ is now $\|\blambda - \blambda_1\|$. 

The vectors $\x,\blambda$ are computed by $\x = \x_k + \delta \x$ and $\blambda = \blambda_k + \delta \blambda$, where the increments are computed via the optimality conditions of~\cref{qp1}, 
\begin{equation}\label{x.l.1}
\Q_k \delta\x + h\, \c_k + \A^\transp \delta \blambda  = \mathbf{0}, \qquad \A \delta \x = {\bf 0}.
\end{equation}

In addition, define the exact prediction computed starting from $(\x^*(t_k), \blambda^*(t_k))$ as $(\x^*_{k+1|k}, \blambda^*_{k+1|k})$ and the increments $\delta \x^* = \x^*_{k+1|k} - \x^*(t_k) = \x^*_{k+1|k} - \x_0$ and $\delta \blambda^* = \blambda^*_{k+1|k} - \blambda^*(t_k)= \blambda^*_{k+1|k} - \blambda_0$, which are computed by [Cf.~\cref{x.l} or equivalently~\cref{qp}] 
\begin{equation}\label{x.l.11}
\Q_0 \delta \x^* + h\, \c_0 + \A^\transp \delta \blambda^* = \mathbf{0}, \qquad \A \delta \x^* = {\bf 0}.
\end{equation} 

The error $\|\x - \x_1\|$ can be upper bounded  as
\begin{multline}\label{er1}
\|\x - \x_1\| \leq \|\delta \x + \x_k - (\delta \x^*+\x_0)\| +\|\x^*_{k+1|k}- \x_1\|  \\ \leq \|\x_k - \x_0\| +  \|\delta \x  - \delta \x^*\| + \Delta_1\Delta_3 h^2, 
\end{multline}
and similarly the error $\|\blambda - \blambda_1\|$ can be upper bounded  as
\begin{equation}\label{er2}
\|\blambda - \blambda_1\| \leq  \|\blambda_k - \blambda_0\| + \|\delta \blambda - \delta \blambda^*\| + \Delta_2\Delta_3 h^2. 
\end{equation}


Consider the solution mapping:
\begin{equation}
s(\p) :=\left\{\y,\bmu\in\im(\A)\left| \begin{array}{c}\Q_k \y + h\, \c_k + \A^\transp \bmu  + \p = \mathbf{0}\\ \A \y=\mathbf{0}\end{array} \right.\right\}.
\end{equation}
The mapping $s(\p)$ is every-where single-valued (due to Proposition~\ref{prop:properties}). In addition, by looking at two different parameters $\p'$ and $\p''$ as done similarly in~\cref{p.pert} and by using Proposition~\ref{prop.perturbation}, then we derive that everywhere (i.e., for all $\p',\p''$), $\|\y(\p') - \y(\p'')\| \leq \Delta_1 \|\p'-\p''\|$ and $\|\bmu(\p') - \bmu(\p'')\| \leq \Delta_2  \|\p'-\p''\|$.

Set $\p' = \mathbf{0}$ and $\p'' = (\Q_0-\Q_k) \delta\x^* + h\, (\c_0-\c_k)$, so that $\y(\p') = \delta\x$, $\bmu(\p') = \delta\blambda$, and $\y(\p'')=\delta\x^*$, and $\bmu(\p'')=\delta\blambda^*$. Then, 
\begin{equation}\label{proofB.dummy6}
\|\p'-\p''\| = \|(\Q_k-\Q_0) \delta\x^* + h\, (\c_k-\c_0)\|.
\end{equation}
We proceed now to bound $\|(\Q_k-\Q_0) \delta\x^* + h\, (\c_k-\c_0)\|$, by using Assumption~\ref{as.4}
\begin{multline}\label{g1_ex}
\|(\Q_k-\Q_0) \delta\x^* + h\, (\c_k-\c_0)\| \leq \\ \leq C_1\|\x_k - \x_0\| \|\delta\x^*\| + h C_2 \|\x_k - \x_0\| .
\end{multline}

The next step is to upper bound $\|\delta\x^*\|$. For this purpose, we use Proposition~\ref{prop.perturbation} on Problem~\cref{x.l.11}. In particular, by~\cref{xx}, one has $\|\delta\x^*\| \leq \Delta_1 C_0 h$, and therefore, 
\begin{align}\label{last}
\|\p'-\p''\|\leq  h (\Delta_1 C_1 C_0 + C_2) \|\x_k - \x_0\|. 
\end{align}
By putting together~\cref{last} with the fact that $\|\y(\p') - \y(\p'')\| \leq \Delta_1 \|\p'-\p''\|$ and $\|\bmu(\p') - \bmu(\p'')\| \leq \Delta_2  \|\p'-\p''\|$ and with~\cref{er1}-\cref{er2}, the bounds~\cref{eq.inside} are proven.
\end{proof}

\subsection{Main algorithm's convergence}

We divide the proof in different steps. Step 1: we bound the prediction error by using Proposition~\ref{prop.intheo.1}; Step 2: we bound the correction error; Step 3: we put the previous steps together and derive the convergence requirements and results.

\textbf{Prediction error.} The distance between the approximate prediction $(\widehat{\x}_{k+1|k},\widehat{\blambda}_{k+1|k})$ and the exact prediction $({\x}_{k+1|k},{\blambda}_{k+1|k})$ can be bounded by using Theorem~\ref{theo.dualascent}. First, notice that, for Theorem~\ref{theo.dualascent} applied to iterations~\cref{dual.ascent.tv}, one has
\begin{subequations}\label{v.result.step1}
\begin{align}
\|\delta{\blambda}_{P} - \delta{\blambda}\| & \leq \varrho_{\textrm{P}}^{P} \|\delta{\blambda}_{0} - \delta{\blambda}\|, \\
\|\delta{\x}_{P} - \delta{\x}\| & \leq \frac{\sigma_{\max}}{m}\varrho_{\textrm{P}}^{P-1}  \|\delta{\blambda}_{0} - \delta{\blambda}\|,
\end{align}
\end{subequations}
or equivalently, by putting $\delta \blambda_0 = {\bf 0}$,
\begin{subequations}\label{v.result}
\begin{align}
\|\widehat{\blambda}_{k+1|k} - {\blambda}_{k+1|k}\| & \leq \varrho_{\textrm{P}}^{P} \|{\blambda}_{k} - {\blambda}_{k+1|k}\|, \\
\|\widehat{\x}_{k+1|k} - {\x}_{k+1|k}\| & \leq \frac{\sigma_{\max}}{m}\varrho_{\textrm{P}}^{P-1}  \|\blambda_{k} - {\blambda}_{k+1|k}\|.
\end{align}
\end{subequations}

By putting together Proposition~\ref{prop.intheo.1}, \cref{v.result}, and~\cref{bounds.xl}, we obtain for the total error after prediction for the dual variable as
\begin{align}
\|\widehat{\blambda}_{k+1|k} - {\blambda}^*(t_{k+1})\| &\leq \|\widehat{\blambda}_{k+1|k} - {\blambda}_{k+1|k}\| +\nonumber \\ & \hspace{1cm}+ \|{\blambda}_{k+1|k} - {\blambda}^*(t_{k+1})\| \leq \nonumber \\
 &\hspace{-2.7cm} \leq \varrho_{\textrm{P}}^{P} \|{\blambda}_{k} - {\blambda}_{k+1|k}\| + \|{\blambda}_{k+1|k} - {\blambda}^*(t_{k+1})\| \nonumber \\
 &\hspace{-2.7cm} \leq  \varrho_{\textrm{P}}^{P}( \|{\blambda}_{k} - {\blambda}^*(t_k)\| + \|{\blambda}^*(t_k) - {\blambda}^*(t_{k+1})\| + \nonumber \\ & \hspace*{-2cm} \|{\blambda}^*(t_{k+1}) - {\blambda}_{k+1|k}\|) +\|{\blambda}_{k+1|k} - {\blambda}^*(t_{k+1})\|\nonumber \\ 
&\hspace{-2.7cm} \leq \varrho_{\textrm{P}}^{P} \|{\blambda}_{k} - {\blambda}^*(t_k)\| + ( \varrho_{\textrm{P}}^{P}+1)\|{\blambda}_{k+1|k} - {\blambda}^*(t_{k+1})\| + \nonumber \\ & \hspace*{1cm} +  \varrho_{\textrm{P}}^{P} \|{\blambda}^*(t_k) - {\blambda}^*(t_{k+1})\| \nonumber \\ 
&\hspace{-2.7cm} \leq  \alpha_2\|{\blambda}_{k} - {\blambda}^*(t_k)\| + h \alpha_1\|\x_{k} - \x^*(t_k)\| + h \alpha_0, \label{final.result}
\end{align}
where we have set $\alpha_2 = 2\varrho_{\textrm{P}}^{P}+1 $, $\alpha_1 = \Delta_2 \Delta_4 (\varrho_\mathrm{P}^{P}+1) $, and $\alpha_0 = \varrho_{\textrm{P}}^{P}(\Delta_2 \Delta_3 h + \Delta_2 C_0) + \Delta_2 \Delta_3 h $.
\vskip5mm

\noindent \textbf{Correction error.} We look now at the correction step, which by using Theorem~\ref{theo.dualascent}, one can derive
\begin{subequations}\label{corr_eq}
\begin{align}
\|{\blambda}_{k+1} - {\blambda}^*(t_{k+1})\| & \leq \varrho_{\textrm{C}}^{C} \|\widehat{\blambda}_{k+1|k} - {\blambda}^*(t_{k+1})\|, \\
\|{\x}_{k+1} - {\x}^*(t_{k+1})\| & \leq \frac{\sigma_{\max}}{m}\varrho_{\textrm{C}}^{C-1}  \|\widehat{\blambda}_{k+1|k} - {\blambda}^*(t_{k+1})\|.
\end{align}
\end{subequations}
with $\varrho_{\textrm{C}} = \max\{|1-\alpha m|,|1 - \alpha L|\}$. And by putting together the result~\cref{final.result} with~\cref{corr_eq}, we obtain the error bounds,
\begin{subequations}\label{rec}
\begin{align}
\|{\blambda}_{k+1} - {\blambda}^*(t_{k+1})\| &\leq \varrho_{\textrm{C}}^{C}(\alpha_2\|{\blambda}_{k} - {\blambda}^*(t_k)\| + \nonumber \\ &\hspace{.5cm}+ h \alpha_1\|\x_{k} - \x^*(t_k)\| + h \alpha_0),\\
\|{\x}_{k+1} - {\x}^*(t_{k+1})\| &\leq \frac{\sigma_{\max}}{m}\varrho_{\textrm{C}}^{C-1}  (\alpha_2\|{\blambda}_{k} - {\blambda}^*(t_k)\| +  \nonumber \\ &\hspace{.5cm}+ h \alpha_1\|\x_{k} - \x^*(t_k)\| + h \alpha_0).
\end{align}
\end{subequations}

\noindent \textbf{Global error and convergence.} Call $a_1 : = \varrho_{\textrm{C}}^{C}\alpha_2$, $a_2 : = h \varrho_{\textrm{C}}^{C}\alpha_1$, $u_1:= h \varrho_{\textrm{C}}^{C} \alpha_0$, and $\gamma := \sigma_{\max}/(m\varrho_{\textrm{C}})$. Define $z_{\blambda,k}: = \|{\blambda}_{k} - {\blambda}^*(t_{k})\|$ and $z_{\x,k}: = \|{\x}_{k} - {\x}^*(t_{k})\|$. Then the error dynamics~\cref{rec} can be written -- in the worst case --  as
\begin{equation}\label{rec1}
\left[\begin{array}{c}z_{\blambda,k+1} \\ z_{\x,k+1} \end{array}\right] = \left[\begin{array}{cc} a_1 & a_2 \\ \gamma a_1 & \gamma a_2  \end{array}\right]\left[\begin{array}{c}z_{\blambda,k} \\ z_{\x,k} \end{array}\right] + \left[\begin{array}{c}u_1 \\ \gamma u_1 \end{array}\right].
\end{equation}
Asymptotic stability of the linear system~\cref{rec1} is achieved iff the eigenvalues of the state transition matrix are inside the unit circle, i.e., iff
\begin{multline}
a_1 + \gamma a_2 < 1 \quad \textrm{i.e.,} \quad \varrho_{\textrm{C}}^{C}(2\varrho_\textrm{P}^{P}+1) +\\+  \sigma_{\max}/(m\varrho_{\textrm{C}}) ( h \varrho_{\textrm{C}}^{C}\Delta_2 \Delta_4 (\varrho_\textrm{P}^{P}+1) ) = \tau(h)<1,  
\end{multline}
that is 
\begin{equation}
h < \frac{m}{\sigma_{\max}}\left[1 - \varrho_{\textrm{C}}^{C}(2\varrho_\textrm{P}^{P}+1)\right]\left[\varrho_{\textrm{C}}^{C-1}\Delta_2 \Delta_4 (\varrho_\textrm{P}^{P}+1) \right]^{-1},
\end{equation}
which is condition~\cref{cond:2} when defining
\begin{equation}
\gamma_1 = \varrho_{\textrm{C}}^{C}(2\varrho_\textrm{P}^{P}+1),
\end{equation}
\begin{multline}
\gamma_2 = \frac{\sigma_{\max}}{m}\left[\varrho_{\textrm{C}}^{C-1}\Delta_2 \Delta_4 (\varrho_\textrm{P}^{P}+1) \right] = \\ \varrho_{\textrm{C}}^{C-1} \frac{\kappa_f \kappa_{\A}^2}{m}\left(\frac{\kappa_f \kappa_{\A}^2 + 1}{m} C_1 C_0 + C_2\right) (\varrho_\textrm{P}^{P}+1) .
\end{multline}
A positive (and therefore implementable) sampling period $h$ exists iff 
\begin{equation}
1 - \varrho_{\textrm{C}}^{C}(2\varrho_\textrm{P}^{P}+1) >0,
\end{equation}
which is condition~\cref{cond:1}, and in this case, 
\begin{multline}\label{rec11}
\left[\!\!\begin{array}{c}z_{\blambda,k} \\ z_{\x,k} \end{array}\!\!\right] \!= \!\left[\!\!\begin{array}{cc} a_1 & a_2 \\ \gamma a_1 & \gamma a_2  \end{array}\!\!\right]^k \left[\!\!\begin{array}{c}z_{\blambda,0} \\ z_{\x,0} \end{array}\!\!\right] \!+\! \sum_{\tau = 0}^{k-1}\left[\!\!\begin{array}{c} (a_1 + \gamma a_2)^\tau \\ \gamma (a_1 + \gamma a_2)^\tau \end{array}\!\!\right] \!u_1 \\ < \infty.
\end{multline}
The asymptotical error is achieved exponentially fast and it is
\begin{subequations}\label{as_error_00}
\begin{align}
\limsup_{k\to\infty} \|{\blambda}_{k} - {\blambda}^*(t_{k})\| & = \frac{u_1}{1-(a_1 + \gamma a_2)} = \nonumber \\ &  \hspace{-2cm} =\frac{\varrho_{\textrm{C}}^{C}[\varrho_{\textrm{P}}^{P}(\Delta_2 \Delta_3 h + \Delta_2 C_0) + \Delta_2 \Delta_3 h]\,h }{1 - \tau(h)},\\
\limsup_{k\to\infty} \|{\x}_{k} - {\x}^*(t_{k})\| & = \frac{\gamma u_1}{1-(a_1 + \gamma a_2)} \nonumber \\ & \hspace{-2cm} = \frac{\sigma_{\max}\varrho_{\textrm{C}}^{C-1} [\varrho_{\textrm{P}}^{P}(\Delta_2 \Delta_3 h + \Delta_2 C_0) + \Delta_2 \Delta_3 h]\,h }{[1 - \tau(h)]m}.
\end{align}
\end{subequations}
Which concludes the proof. \qed

\section{Proof of Claim~\ref{claim:dist}}\label{appendix:cinque}

\begin{proof}
To justify the claim, we analyze all the steps of Algorithm~\ref{alg.approx}. First, the prediction step is based on the iterations~\cref{dual.ascent.tv}. Let $\y_{i,k}$ and $\delta \y_{i,k}$ be the local variables $\y_i$ and $\delta \y_i$ at iteration $k$; let $\delta \blambda_{i,j,k}$ be the dual variable associated with link $(i,j)\in E$ at iteration $k$, then~\cref{dual.ascent.tv} can be rewritten as

\begin{enumerate}
\item For all $i\in V$ do
\begin{subequations}\label{dual.ascent.tv.dist}
\begin{multline}
\delta \y_{i,p+1} = \argmin_{\delta \y_i \in \mathbb{R}^n}\Big\{\frac{1}{2}\delta\y_i^\transp\nabla_{\y_i\y_i} f_i(\y_{i,k}; t_k) \delta\y_i + \\ + \nabla_{t\y_i}f(\y_{i,k}; t_k)^\transp \delta\y_i+ \\ \sum_{(i,j)\in E, i\leq j}\!\!\delta \blambda_{i,j,p}^\transp \delta \y_{i} - \sum_{(i,j)\in E, i> j}\!\! \delta \blambda_{i,j,p}^\transp \delta \y_{i} \Big\},
\end{multline}
\item Communicate $\delta \y_{i,p+1}$ with neighbors;
\item For all $i \in V$ do
\begin{align}
\delta \blambda_{i,j,p+1}& = \delta \blambda_{i,j, p} + \beta \mathbb{I}_{i\leq j}(\delta \y_{i,p+1}-\delta\y_{j,p+1}),
\end{align}
\end{subequations}
where $\mathbb{I}_{i\leq j}$ is $1$ if $i\leq j$, and $-1$ otherwise. 
\end{enumerate}
This justifies the fact that the prediction step can be implemented in a distributed fashion with synchronous communication (Assumption~\ref{as.synch}). Each node maintains local copies $\delta \y_{i,k}, \delta \blambda_{i,j,k}$ which converge to the primal-dual optimizers of the prediction step. 

Second, we analyze the correction step, which is based on the iterations~\cref{dual.ascent.tv.correction}. It is easy to see that also~\cref{dual.ascent.tv.correction} can be written in a similar fashion as~\cref{dual.ascent.tv.dist}, thereby allowing distributed computation of the correction direction with synchronous communication (Assumption~\ref{as.synch}).

Provided now that the switching between prediction and correction step is synchronized (Assumption~\ref{as.synch}) then Algorithm~\ref{alg.approx} can be implemented in a distributed fashion. \end{proof}

\section{Asymptotical error bounds}\label{appendix:sei}

We prove here both~\eqref{err.cec} and \eqref{err.tc}. For~\eqref{err.cec}, by similar arguments as the one of the proof of Theorem~\ref{theo.dualascent.tv}, we can write, 
\begin{align}
\|\blambda_{k+1} - \blambda^*(t_{k+1})\| &\leq \varrho_{\textrm{C}}^C (\|\widetilde{\blambda_{k}} - \blambda^*(t_{k})\| + K)  \nonumber \\ & \leq \varrho_{\textrm{C}}^C ( \varrho_{\textrm{C}}^{C'}\|{\blambda_{k}} - \blambda^*(t_{k})\| + K)  \nonumber \\
\|\x_{k+1} - \x^*(t_{k+1})\| &\leq \frac{\sigma_{\max}}{m} \varrho_{\textrm{C}}^{C-1}  (\varrho_{\textrm{C}}^{C'}\|{\blambda_{k}} - \blambda^*(t_{k})\| + K),
\end{align}
and therefore,
\begin{multline}
\mathrm{Err}_{\mathrm{C+EC}} = \limsup_{k\to\infty} \|\x_{k} - \x^*(t_{k})\| = \\\frac{\sigma_{\max}}{m} \varrho_{\textrm{C}}^{C-1}  \left(\frac{\varrho_{\textrm{C}}^{C+C'} K}{1-\varrho_{\textrm{C}}^C} + K\right),
\end{multline}
from which~\eqref{err.cec}. 

As for~\eqref{err.tc}, 
\begin{align}
\|\blambda_{k+1} - \blambda^*(t_{k+1})\| &\leq \varrho_{\textrm{C}}^{C''} (\|{\blambda_{k}} - \blambda^*(t_{k})\| + K)  \nonumber \\ 
\|\x_{k+1} - \x^*(t_{k+1})\| &\leq \frac{\sigma_{\max}}{m} \varrho_{\textrm{C}}^{C''-1}  (\|{\blambda_{k}} - \blambda^*(t_{k})\| + K),
\end{align}
and therefore,
\begin{multline}
\mathrm{Err}_{\mathrm{TC}} = \limsup_{k\to\infty} \|\x_{k} - \x^*(t_{k})\| = \\\frac{\sigma_{\max}}{m} \varrho_{\textrm{C}}^{C''-1}  \left(\frac{\varrho_{\textrm{C}}^{C''} K}{1-\varrho_{\textrm{C}}^{C''}} + K\right),
\end{multline}
from which~\eqref{err.tc}.

\bibliographystyle{ieeetr}
\bibliography{../../../../../PaperCollection00}

\end{document}